\input amstex
\def\sk1{\vskip 10pt}
\def\newline{\hfil\break} 
\def\cl{\centerline}
\magnification\magstep1
\mathsurround=1pt
\baselineskip=20pt plus1pt

\def\eset{\emptyset}

\def\spitem{\par\hangone\textputone}
\def\hangone{\hangindent 3\parindent}
\def\textputone#1{\noindent
	 \hbox to 3\parindent{\hss#1\enspace}\ignorespaces}

\def\ds{\displaystyle}

\def\ni{\noindent}

\cl{A Proof Of The Hilbert-Smith Conjecture}

\cl{by}

\cl{Louis F\. McAuley}
\cl{Dedicated to the memory of Deane Montgomery}
\vskip 30pt

\cl{Abstract}

The Hilbert-Smith Conjecture states that if $G$ is a locally compact group
which acts effectively on a connected manifold as a topological
transformation group, then $G$ is a Lie group.  A rather straightforward proof of
this conjecture is given.  The motivation is work of Cernavskii
(``Finite-to-one mappings of manifolds'', {\it Trans\. of Math}\. Sk\. 65
(107), 1964.)  His work is generalized to the orbit map of an
effective action of a $p$-adic group on compact connected
$n$-manifolds with the aid of some new ideas.  There is no attempt to
use Smith Theory even though there may be similarities.
\vskip 20pt

\ni 1.  {\it Introduction}.  

In 1900, Hilbert proposed twenty-three problems [8].  For an excellent
discussion concerning these problems, see the {\it Proceedings of
Symposia In Pure Mathematics} concerning ``Mathematical Developments
Arising From Hilbert Problems'' [3].  The abstract by C.T. Yang [22]
gives a review of Hilbert's Fifth Problem ``{\it How is Lie's concept of
continuous groups of transformations of manifolds approachable in our
investigation without the assumption of differentiability}?''  Work of von
Neumann [40] in 1933 showed that differentiability is not completely
dispensable.  This with results of Pontryagin [35] in 1939 suggested the
specialized version of Hilbert's problem:  {\it If} $G$ {\it is a
topological group and a topological manifold, then is} $G$ {\it
topologically isomorphic to a Lie group}?  This is generally regarded as
Hilbert's Fifth Problem.  The first partial result was given by Brouwer [26]
in 1909-1910 for locally euclidean groups of dimension $\leq 2$.  The best
known partial results were given for compact locally euclidean groups and
for commutative locally euclidean groups by von Neumann [40] and Pontryagin
[35], respectively.

In 1952, the work of Gleason [5] and Montgomery-Zippin [13] proved:  {\it
Every locally euclidean group is a Lie group}.  This solved Hilbert's Fifth
Problem.  

A more general version of Hilbert's Fifth Problem is the following:

{\it If} $G$ {\it is a locally compact group which acts effectively on a connected
manifold as a topological transformation group, then is} $G$ {\it a Lie
group}?  The Hilbert-Smith Conjecture states that the answer is yes.  

Papers of Montgomery [34] in 1945 and Bochner-Montgomery [1] in 1946
established the partial result:  {\it Let} $G$ {\it be a locally compact
group which acts effectively on a differentiable manifold} $M$ {\it such
that for any} $g \in G$, $x \mapsto gx$ {\it is a differentiable
transformation of} $M$.  {\it Then} $G$ {\it is a Lie group and} $(G,M)$
{\it is a differentiable transformation group}.  Another partial result was
given by a theorem of Yamabe [43] and a Theorem of Newman [15] as follows:
{\it If} $G$ {\it is a compact group which acts effectively on a manifold
and every element of} $G$ {\it is of finite order, then} $G$ {\it is a
finite group}.  

It has been shown [14] that an affirmative answer to the generalized version
of Hilbert's Fifth Problem
is equivalent to a negative answer to the following:  {\it Does there exist
an effective action of a} $p$-{\it adic group on a manifold}?

It is proved here that the answer to this question is No.
Thus, the {\it
Hilbert-Smith Conjecture} is true, 
i.e., {\it A locally compact group acting
effectively on a connected $n$-manifold must be a Lie group}. 

A brief review of some of the consequences of efforts to solve this problem
is given below.  There are examples in the literature of {\it effective}
actions of an infinite compact 0-dimensional topological group $G$ ({\it
each} $g \in G$-$\{$identity$\}$ 
{\it moves some point}) on locally connected continua.  The
classic example of Kolmogoroff [29] in 1937, is one where $G$ operates
effectively but not {\it strongly effectively} [24] 
on a 1-dimensional locally connected continuum (Peano continuum) such
that the orbit space is 2-dimensional.  In 1957, R.D\. Anderson [24] proved
that {\it any} compact 0-dimensional topological group $G$ can act strongly
effectively as a transformation group on the (Menger) universal
1-dimensional curve $M$ such that either (1) the orbit space is homeomorphic
to $M$ or (2) the orbit space is homeomorphic to a regular curve.  

In 1960, C.T\. Yang [44] proved that if a $p$-adic group, $A_p$,
acts effectively as
a transformation group on $X$ (a locally compact Hausdorff space of homology
dimension not greater than $n$), then the homology dimension of the orbit
space $X/A_p$ is not greater than $n + 3$.  If $X$ is an $n$-manifold, then
the homology dimension of $X/A_p$ is $n + 2$.  If $A_p$ acts strongly
effectively (freely) on an $n$-manifold $X$, then the dimension of $X/A_p$ is
either $n + 2$ or infinity.  At about the same time (1961), Bredon, Raymond, and
Williams [25] proved the same results using different methods.  There are,
of course, actions by $p$-adic groups on $p$-adic solenoids and actions by
$p$-adic solenoids on certain spaces.  See [25] for some of these results.

In 1961, Frank Raymond published the results of his study of the orbit space
$M/A_p$ {\it assuming} an effective action by $A_p$ (as a transformation
group) on an $n$-manifold $M$.  Later (1967), Raymond [38] published work on
two problems in the theory of generalized manifolds which are related to
the (generalized) Hilbert Fifth Problem.

In 1963, Raymond and Williams [39] gave examples of compact metric spaces
$X^n$ of dimension $n$ and an action by a $p$-adic group, $A_p$, on $X^n$
such that dim $X^n/A_p = n + 2$.  Work related to and used in [39] is the
paper [41] by Williams.  In [41], Williams answers a question of Anderson
[24; p\. 799] by giving a free action by a compact 0-dimensional group $G$
on a 1-dimensional Peano continuum $P$ with dim $P/G = 2$. 

In 1976, I described [32; 33] what I called $p$-adic polyhedra which admit
periodic homeomorphisms of period $p$.  Proper inverse systems
$\{P_i,\phi_i\}$ of $p$-adic $n_i$-polyhedra have the property that the
inverse limit $X = \ds{\varprojlim} P_i$ admits a free action by a
$p$-adic group.  

In 1980, one of my students, 
Alan J\. Coppola [28] generalized results of C.T\. Yang [44] which
involve homologically analyzing $p$-adic actions.  Coppola formalized these
so that homological calculations could be done in a more algorithmic manner.
He defined a $p$-adic transfer homomorphism and used it to produce all of
the relevant Smith-Yang exact sequences which are used to homologically
analyze $Z_{p^r}$-actions on compact metric spaces.  Coppola studied
$p$-adic actions on homologically uncomplicated spaces.  In particular, he
proved that if $X$ is a compact metric $A_p$-space of homological dimension
no greater than $n$ and $X$ is homologically locally connected, then the
$(n+3)$-homology of any closed subset $A \subset X/A_p$ vanishes.

In 1983, Robinson and I proved Newman's Theorem for finite-to-one open and
closed mappings on manifolds [10].  We formalized Newman's Property (and
variations) and studied this property for discrete open and closed mappings
on generalized continua in 1984 [11].  

In 1985, H-T Ku, M-C Ku, and Larry Mann investigated in [30] the connections
between Newman's Theorem involving the size of orbits of group actions on
manifolds and the Hilbert-Smith Conjecture.  They establish Newman's Theorem
(Newman's Property [11]) for actions of compact {\it connected} non-Lie
groups such as the $p$-adic solenoid.  

In 1997, D.\ Repov{\v s} and E.V.\ {\v S}{\v c}epin [51] gave a proof of the Hilbert-Smith
Conjecture for actions by Lipschitz maps.  See also related work by Shchepin [52].  In
the same year, Iozke Maleshick [53] proved the Hilbert-Smith conjecture for H{\"o}lder
actions.

In 1999, Gaven J. Martin [54] announced a proof of The Hilbert-Smith Conjecture for
quasiconformal actions on Riemannian manifolds and related spaces.

The crucial idea that works here is M.H.A. Newman's idea used in his proof that
for a given compact connected $n$-manifold $M$, 
there is an
$\epsilon > 0$ such that if $h$ is any periodic homeomorphism of period
$p$, a prime $> 1$, of $M$ onto itself, then
there is some $x \in M$ such that the orbit of $x$, 
$\{x,h(x),\hdots,h^{p-1}(x)\}$, has diameter $\geq \epsilon$. 
It is well known that the
collection of orbits under the action of a transformation group $G$ on a
compact Hausdorff space $X$ is a {\it continuous decomposition} of $X$.  

The works [20; 21] of David Wilson and John Walsh [18] show that there exist
continuous decompositions of $n$-manifolds $M^n$, $n \geq 3$, into Cantor
sets.  This paper shows that  
such decompositions can not be equivalent to those induced 
by any action of a $p$-adic transformation group $A_p$ acting on $M^n$.

I owe a special debt of gratitude to Patricia Tulley McAuley who has
been extremely helpful in reading drafts of numerous attempts to solve this
problem and who has
provided helpful insights with regard to {\v C}ech homology.  Also, I am
indebted to the work of Cernavskii [27].  
\sk1

\cl{\it OUTLINE OF A PROOF}

It is well know that if a locally compact group $G$ acts effectively on a
connected $n$-manifold $M$ and $G$ is not a Lie group, then there is a subgroup
$H$ of $G$ isomorphic to a $p$-adic group $A_p$ which acts effectively on $M$.
Thus, the Hilbert-Smith Conjecture can be established by proving that there is no
effective action by a $p$-adic group $A_p$ on a connected $n$-manifold $M$.  The
conjecture is proved by the following theorem.  As seen later, there is
no loss of generality, in assuming that $M$ is compact, orientable, and without
boundary.   The proof given below can be adopted to the situation where $M$ is a
{\it locally compact}, orientable, and without boundary by replacing finite open
coverings of $M$ and $M/A_p$ with locally finite open coverings with the same
properties since the orbit map $\phi : M\to M/A_p$ is open, closed, and proper.

\ni Theorem. If $L(M,p)$ is the class of all orbit mappings $\phi : M \to
M/A_p$ where $A_p$ acts effectively on a compact, connected, and orientable
$n$-manifold $M$, then $M$ has Newman's Property w.r.t. $L(M,p)$.  

It is well known that $M$ does not have Newman's Property w.r.t. $L(M,p)$.
Hence, $L(M,p) = \eset$ and the Hilbert-Smith Conjecture is true. 

\ni Definition.
{\rm An $n$-manifold $(M,d)$ is said to have Newman's
Property w.r.t\. the class $L(M,p)$ (as stated above) iff 
there is $\epsilon > 0$ such that for any $\phi \in L(M;p)$, there is 
some $x\in
M$ such that diam $\phi^{-1}\phi(x) \geq \epsilon$ using the metric $d$ on $M$.}

\ni Lemma 2. (A consequence of a Theorem of Floyd [4].)
Suppose that $M$ is a compact connected orientable
$n$-manifold.  There is a finite open covering
$W_1$ of $M$ such that (1) order $W_1 = n+1$ and
(2) there is a finite open refinement $W_2$ of $W_1$ which covers $M$ such
that if $W$ is any finite open covering of $M$ refining $W_2$,
then $\pi_{W_1} : {\check H}_n(M) \to H_n(W_1)$ maps ${\check H}_n(M)$
isomorphically onto the image of the projection $\pi_{WW_1} : H_n(W) \to
H_n(W_1)$.  

[Here, if $U$ is either a finite or locally finite
open covering of $M$, then $H_n(U)$ is the
$n^{\text{th}}$ simplicial homology group of the nerve $N(U)$ of $U$.  The
coefficient group is always $Z_p$ and ${\check H}_n(M)$ denotes the
$n^{\text{th}}$ {\v C}ech homology of $M$.]    

Now, choose $U = W_1$, and a finite open covering $W$ of $M$ 
which star refines $W_2$ where $W_1$ and $W_2$ satisfy Lemma 2.

Let $\epsilon$ be the Lebesque number of
$W_2$.  Choose $\phi \in L(M,p)$ such that diam $\phi^{-1}\phi(x) < \epsilon$ for each
$x\in M$.   Construct the special
coverings $\{V^m\}$ and the special refinements $\{U^m\}$ as in 
Lemma 5 below where $V^1$ star refines $W_2$ such that order $U^m = n+1$ and with 
projections $\alpha_m$, $\beta_m$, and $\pi_m$ yielding the
following commutative diagram: 
$$
\matrix
\gets &H_n(V^m) &\overset {\beta^*_m} \to \longleftarrow &H_n(U^m)
&\overset {\alpha^*_m} \to \longleftarrow &H_n(V^{m+1}) &\overset
{\beta^*_{m+1}} \to \longleftarrow &H_n(U^{m+1}) &\gets & \\
& & & & & & & & & \\
 &{\nu_m} \uparrow &{\beta^*_m}\swarrow & &\nwarrow
{\alpha^*_m} &\nu_{m+1} \uparrow &\beta^*_{m+1}
\swarrow & &\nwarrow \alpha^*_{m+1} & \\
& & & & & & & & & \\
\gets &H_n(V^m_n) & &\overset {\pi^*_m} \to \longleftarrow & &H_n(V^{m+1}_n)
& &\overset{\pi^*_{m+1}} \to \longleftarrow &H_n(V^{m+2}_n)
&\gets
\endmatrix
$$
Here $\nu_m$ is the natural map of an $n$-cycle in $H_n(V^m_n)$, the
$n^{\text{th}}$ simplicial homology group of the $n$-skeleton of the
nerve $N(V^m)$ of $V^m$, into its homology
class in $H_n(V^m)$.  The other maps are those induced by the projections
$\alpha_m$, $\beta_m$ and $\pi_m$.  The upper sequence, of course, yields
the {\v C}ech homology group ${\check H}_n(M)$ as its inverse limit.
Furthermore, it can be easily shown, using the diagram, that ${\check
H}_n(M)$ is isomorphic to the inverse limit $G = \varprojlim H_n(V^m_n)$, of
the lower sequence.  Specifically, $\gamma : {\check H}_n(M) \to G$ defined
by $\gamma(\Delta) = \{\beta^*_m(\pi_{U_m}(\Delta))\}$ is an isomorphism of
${\check H}_n(M)$ onto $G$.  We shall use the isomorphism in what follows
and for convenience we shall let $\gamma(\Delta) = \{z^n_m(\Delta)\}$, i.e.
$z^n_m(\Delta) = \beta^*_m(\pi_{U^m}(\Delta)) \in H_n(V^m_n)$.  
The group $\varprojlim H_n(V^m_n)$ is used because the operator
$\sigma$ (introduced later) is applied to actual $n$-cycles rather
than to elements of a homology class.  This is the reason for the
sequence $\{U^m\}$.  Note that there is no attempt to use Smith
Theory.

An operator $\sigma_m$ is defined on the $n$ chains of
$N(V^{m+1})$ for each $m$.  
The operator $\sigma_m$ maps $n$-cycles to $n$-cycles and commutes with the
projections $\pi^*_m : H_n(V^{m+1}_n) \to H_n(V^m_n)$ and, hence, induces 
an automorphism on $\pi_{V^{m+1}}({\check H}_n(M))
\subset H_n(V^{m+1}_n)$.  See Lemmas 6 and 7.  

It is shown in Lemma 10 that int $F_\phi = \eset$ where $F_\phi = \{x\mid \phi^{-1}
\phi(x) = x\}$.  

Distinguished families of $n$-simplices in $N(V^m)$ are defined.  
Now, let $z_m = z_m^n(\Delta)$ where $\Delta$ is the generator of $G\cong Z_p$.  
For each $n$-simplex $\delta^n$ in
$z_m(\Delta)$, there is a unique distinguished family $S^m_j$ of $n$-simplices in 
$N(V^m)$ which contains $\delta^n$.  If $C_j$ is the collection of all 
$n$-simplices in $z_m(\Delta)$ which are in $S^m_j$, then the sum of the 
coefficients  of those members of $C_j$ (as they appear in $z_m(\Delta))$ is 0 mod
$p$.  Take the projection $\pi_{V^mU}$ from
$V^m$ to $U = W_1$. Hence, $\pi_{V^mU}$ has the property that all members of a 
distinguished
family $S^m_j$ of $n$-simplices 
in $N(V^m)$ project to the
same simplex in $N(W_1)$.   Thus, the projection of those
members of $z_m(\Delta)$ which are in $S^m_j$ 
project to the same simplex $\delta_j$ in $N(W_1)$
and the coefficient of $\delta_j$ is 0 mod $p$.   
Thus, $\pi_{V^mU} : H_n(V^m_n) \to H_n(U)$
takes the nontrivial $n$-cycle $z_m(\Delta)$ to the $0-n$ cycle mod $p$.  This
violates the conclusion of Lemma 2.  Thus, $M$ has Newman's Property w.r.t.\ the
class $L(M,p)$.  Hence, $\epsilon$ is a Newman's number and the Theorem is
proved.  

It is well known that 
if $A_p$ acts effectively on a compact connected $n$-manifold
$M$, then given any $\epsilon > 0$, there is an effective action of $A_p$
on $M$ such that diam $\phi^{-1}\phi(x) < \epsilon$ for each $x\in M$.  That
is, $M$ fails to have Newman's property w.r.t. the class $L(M,p)$.  It follows that
$A_p$ can not act effectively on a compact connected $n$-manifold $M$.  {\it
Consequently, the Hilbert-Smith Conjecture is true}.  

Details of the proof follow.
\vskip 30pt

\ni 2.  {\it Some properties of the orbit mapping of an effective action by
$A_p$ on a compact connected orientable $n$-manifold $M$.}

Suppose that $\phi$ is the orbit mapping of a $p$-adic group $A_p$ acting
as a transformation group on an orientable $n$-manifold $M^n = M$ where $p$
is a prime larger than $1$.  By [12; 21], there is a sequence $A_p = H_0
\supset H_1 \supset H_2 \supset \hdots$ of open (and closed) subgroups of
$A_p$ which closes down on the identity $e$ of $A_p$ such that when $j > i$,
$H_i/H_j$ is a cyclic group of order $p^{j-i}$.  Let $h_{ij} : A_p/H_j \to
A_p/H_i$ and $h_i : A_p \to A_p/H_i$ be homomorphisms induced by the
inclusion homomorphisms (quotient homomorphisms) on $A_p$ and $A_p/H_j$.
Then $\{ A_p/H_i;h_{ij}\}$ is an inverse system and $\{h_i\}$
gives an isomorphism of $A_p$ onto $\varprojlim\ A_p/H_i$.  
Now, let $a \in A_p - H_i$.
For each natural number $i$, let $a_i$ be the coset $aH_i$ in $A_p/H_i$.
Then $a_i$ is a periodic homeomorphism of $M/H_i$ onto $M/H_i$ with $a^q_i$
being the identity mapping where $q = p^i$ is the period of $a_i$.
Consequently, $H_i$ acts as a transformation group on $M$ and $A_p/H_i$ acts
as a cyclic transformation group on $M/H_i$.  

As above, let $\{H_i\}$ be a sequence of open (and closed) subgroups of
$A_p$ such that (a) $H_i \supset H_{i+1}$ for each $i$, (b) if $j \geq i$,
then $H_i/H_j$ is a cyclic group of order $p^{j-i}$, and (c) $A_p/H_i$ is a
cyclic group of order $p^i$.  Since $A_p$ acts effectively on $M$ (a compact
connected $n$-manifold), the cyclic group $A_p/H_i$ acts
effectively on $M/H_i$ with orbit space $M/A_p$.  Let $\pi_{ij} : M/H_j \to
M/H_i$ where $j > i$ and $\pi_i : M\to M/H_i$ be maps induced by the
identity map of $M$.  Also, $\pi_i = \pi_{ji} \circ \pi_j$.  Thus, $\{M/H_i
: \pi_{ij}\}$ is an inverse system and $\{\pi_i\}$ gives a homeomorphism of
$M$ onto the inverse limit $M/H_i$ ($\varprojlim M/H_i)$.  Notice that
$H_i/H_j$, a cyclic group of order $p^{j-i}$, acts on $M/H_j$ with orbit
space $M/H_i$ [cf\. 23; p\. 211] where $\pi_i = \pi_{ji} \circ \pi_j$ and
$\pi_{ji} : M/H_j \to M/H_i$ is the orbit map of the action of $H_i/H_j$ on
$M/H_j$.  
Notice that $M/H_j$ is the orbit space of the action of $H_j$ on $M$ where $H_j$ is
isomorphic to $A_p$.  It follows that if $\epsilon > 0$, then there is a natural number $j$
such that $H_j \cong A_p$ acts effectively on $M$ such that if $\phi_j : M\to M/H_j$ is the
orbit map of the action, then diam $\phi^{-1}\phi(x) < \epsilon$ for each $x\in M$.
Observe that if $H$ is a non empty open and closed subgroup of $A_p$,
then for some $i$, $H = H_i$.  

Observe that if $A_p$ acts effectively as a transformation group on $M$
where $M$ is an orientable connected metric $n$-manifold, then some orbit is infinite
and not discrete.  This follows from a Theorem of Cernavskii [27] that if
$\phi : M\to Y$ (Hausdorff) is a discrete open and closed (continuous)
mapping, then there is a natural number $k$ such that cardinality of
$\phi^{-1}\phi(x) \leq k$ for each $x \in M$ (bounded multiplicity).
Furthermore, the union of the orbits of maximal cardinality is a dense open
set $W$ in $M$.  The stability group of all the points in $W$ is a certain
$H_v$ and, thus, $H_v$ acts as the identity on $M$.
Consequently, the action of $A_p$ is not
effective contrary to the hypothesis.  
Thus, the
orbit mapping $\phi : M \to M/A_p$ where $A_p$ acts effectively on $M$ is
not a discrete open and closed mapping.  Hence, some orbit is
infinite and not discrete.  

The following lemma is crucial to defining certain coverings of $M$ with
distinguished families of open sets.  

\ni Lemma 1.  Suppose that $\phi$ is the orbit mapping $\phi : M\to M/A_p$
where $A_p$ acts effectively on a compact orientable connected
$n$-manifold $M$.  For each $z\in M/A_p -
\phi(F_\phi)$ and $\epsilon > 0$ where $F_\phi = \{x\mid \phi^{-1}\phi(x) = x\}$, 
there is a connected set 
$U$ such that (1) diam $U < \epsilon$, (2) $z\in U$, and (3)
$\phi^{-1}(U) = \{U_1,U_2,\cdots,U_{p^s}\}$ where $s$ is a natural number
such that (a) $U_i$ is a
component of $\phi^{-1}(U)$ for each $i$, $1\leq i \leq p^s$, (b) $\bar U_i \cap
\bar U_j = \eset$ for $i\neq j$, (c) $\phi(U_i) = U$ for each $i$, and (d) $U_1$
is homeomorphic to $U_j$ for each $j$, $1 < j\leq p^s$ (by maps compatible
with the projections $\phi \mid U_j$). 
[The homeomorphism taking $U_1$ to $U_j$ is a power
of a fixed element $g\in A_p - H_1$ and is used in Lemma 6.] 

\ni Proof.  Since $z\in M/A_p -\phi(F_\phi)$, $\phi^{-1}(z)$ is non degenerate.
For $\epsilon > 0$, there is a connected open set $U$ such that (1) diam $U <
\epsilon$, (2) $z\in U$, and (3) $\phi^{-1}(U)$ consists of a finite number
(larger than one) of components $U_1,U_2,\cdots,U_m$ such that $\bar U_i\cap \bar
U_j = \eset$ for $i\neq j$, and $\phi(U_i) = U$ for each $i$.  This follows by 
Whyburn's Theory of open mappings
[19, pp\. 78-80]. 

For each $U_j$, a component of $\phi^{-1}(U)$, there is an open and closed
subgroup $G_j$ of $A_p$ which is the largest subgroup of $A_p$ which leaves
$U_j$ invariant and the map induced by $\phi$ maps $U_j/G_j$ onto $U$.
Since $G_j$ is a normal subgroup of $A_p$, $G_i = G_j$ for each $i$ and $j$.
Furthermore, $A_p/G_j$ is a cyclic group of order $p^s$ where $s$ is a natural
number.  There are $p^s$
pairwise distinct components of $\phi^{-1}(U)$.  (See [36:  Lemma 2]).  It
follows that $G_1 = H_i$ for some $i$ where $\{H_i\}$ is the sequence of open and
closed subgroups of $A_p$ which closes down on the identity $e\in A_p$ (mentioned
above) and $s = i$.  Let $a\in A_p - H_i$ such that $aH_i$ generates
the cyclic group $A_p/H_i$.  
For each natural number $i$,
let $a_i$ be the coset $aH_i$ in $A_p/H_i$.  Thus, $a_i$ is a periodic
homeomorphism of $M/H_i$ onto $M/H_i$ with $a^q_i = e$ where $q = p^i$ is
the period of $a_i$.  

Let $f_i : M\to M/H_i$ be the orbit map of the action of $H_i$ on $M$ and
$g_i : M/H_i \to M/A_p$ be the orbit map of the action of the cyclic group
$A_p/H_i$ on $M/H_i$.  That is, $\phi = g_if_i$.  There are $p^i$ cosets
$\big\{v_mH_i\big\}^{p^i}_{m=1}$ where $v_1 = e$ (the identity) such that
for each $x\in U_1$, $\phi^{-1}\phi(x) = \ds{\bigcup^{p^i}_{m=1}}v_mH_i(x)$
where $v_mH_i(x) = \{h(x)\mid h\in v_mH_i\}$, (2) $v_mH_i(x) \in U_m/H_i$, (3)
$v_m$ is an orientation preserving homeomorphism of $M$ onto $M$, and (4) if
$A_p/H_i = \{a_i,a^1_i,a^2_i,\cdots,a^{p^i-1}_i\}$ (a cyclic group), then
there are elements $k_1,k_2,\cdots,k_{p^i}$ where $k_1~=~e$ such that (a)
$k_m(\pi_i(x)) = \pi_i(v_mH_i(x))$ where $\pi_i$ maps $\phi^{-1}(U)$ onto
$\phi^{-1}(U)/H_i$ and (b) $k_m$ maps $U_1/H_i$ homeomorphically onto
$U_m/H_i$ with $k_m = v_mH_i \in A_p/H_i$ which is a homeomorphism of
$M/H_i$ onto $M/H_i$.  Thus, $v_m(x) \in U_m$.  Let $z\in U_m$.  Hence,
$\pi_i(z) \in U_m/H_i$ and $(v_mH_i)^{-1}(\pi_i(z)) = v^{-1}_mH_i(\pi_i(z)) \in
U_1/H_i$.  Consequently, $v^{-1}_mH_i(\pi_i(z)) = \pi_i(v^{-1}_m(z)) \in
U_1/H_i$ which implies that $v^{-1}_m(z) \in U_1$.  Finally,
$v_m(v^{-1}_m(z)) = z$ and $v_m$ maps $U_1$ homeomorphically onto $U_m$. 

Lemma 1 is proved.
\vskip 30pt

\ni 3.  {\it Special coverings, distinguished families, and distinguished
subfamilies}. 

Let $L(M,p) = \{\phi \mid \phi$ is the orbit mapping of an effective action
of a $p$-adic group $A_p$ ($p$ a prime with $p > 1$) on a compact connected
metric orientable $n$-manifold without boundary, 
$\phi : M\to M/A_p\}$.  For each $\phi \in
L(M,p)$, let $F_\phi = \{x\mid x\in M$ and $\phi^{-1}\phi(x) = x\}$, the
fixed set of the action $A_p$.  It would simplify the proof of lemmas which
follow to know that $M$ is triangulable.  Without this knowledge, a theorem
of E.E\. Floyd is used.

\ni {\bf Notation.}  Throughout this paper, ${\check H}_n(X)$ will denote
the $n^{\text{th}}$ {\v C}ech homology group of $X$ with coefficients in
$Z_p$, the integers mod $p$, $p$ a fixed prime larger than 1.  Also,
$H_n(K)$ will denote the $n^{\text{th}}$ simplicial homology of a finite
simplicial complex $K$, with coefficients in $Z_p$.  If $U$ is a finite open
covering of a space $X$, then $N(U)$ denotes the nerve of $U$, $H_n(U)$ is the
$n^{\text{th}}$ simplicial homology group of $N(U)$, and $\pi_U$ the usual
projection homomorphism $\pi_U : {\check H}_n(X) \to H_n(U)$.  

\ni{\bf Definition.}  If $f$ is a mapping of $M$ onto $Y$, then an open
covering $U$ of $M$ is said to be {\it saturated} (more precisely, saturated
w.r.t\. $f$) iff for each $u\in U$, $f^{-1}f(u) = u$.  That is, $u$ is an
open inverse set.

The next lemma follows.

\ni Lemma 2.  Suppose that $M$ is a compact connected orientable
metric $n$-manifold.  There is a finite open covering
$W_1$ of $M$ such that (1) order $W_1 = n+1$ and
(2) there is a finite open refinement $W_2$ of $W_1$ which covers $M$ such
that if $W$ is any finite open covering of $M$ refining $W_2$,
then $\pi_{W_1} : {\check H}_n(M) \to H_n(W_1)$ maps ${\check H}_n(M)$
isomorphically onto the image of the projection $\pi_{WW_1} : H_n(W) \to
H_n(W_1)$.  

\ni Proof.  Adapt Theorem (3.3) of [4] to the situation here and use (2.5)
of [4].  

If $M$ is triangulable, then there is a sufficiently fine
triangulation $T$ such that if $U$ consists of the open stars of the vertices
of $T$, then $\pi_U : {\check H}_n(M) \to H_n(U)$ is an isomorphism onto
(where, of course, ${\check H}_n(M) \cong Z_p)$.  

{\bf NOTE:}  The reader can assume, for convenience, that $M$ is
triangulable and that $W_1$ is the collection of open stars of a
sufficiently fine triangulation $T$ such that the $n^{\text{th}}$ simplicial
homology, $H_n(W_1)$, with coefficients in $Z_p$, of the nerve, $N(W_1)$, of
$W_1$ is $Z_p$.  

{\it Standing Hypothesis}:  In the following, $M$ is a compact, connected, and
orientable metric $n$-manifold.  Also, $L(M,p)$ is as defined above.  The finite
open coverings $W_1$ and $W_2$ which satisfy Lemma 2 will be used in certain
lemmas and constructions which follow.  Suppose also that $Y = M/A_p$ has a
countable basis $Q = \big\{B_i\big\}^\infty_{i=1}$ such that (a) for each $i$,
$B_i$ is connected and uniformly locally connected and (b) if $H$ is any
subcollection of $Q$ and $\ds{\bigcap_{h\in H}}h\neq \eset$, then
$\ds{\bigcap_{h\in H}}h$ is connected and uniformly
locally connected (a consequence of a
theorem due to Bing and Floyd [50]).  

\ni Lemma 3.  Suppose that $\phi \in L(M,p)$ and $F = F_\phi$.
Then there is a finite open irreducible covering $W_F$ which covers $F$
and refines $W_2$ such that (1) $W_F$ star refines $W_2$ (2) order $W_F \leq
n+1$, (3) for $w\in W_F$, $\phi^{-1}\phi(w) = w$ ($w$ is an open inverse set or
$W_F$ is saturated), (4) if $BdF = B = F -
\text{interior }F$ and $W_B = \{w\mid w \in W_F$ and $w\cap B \neq \eset\}$,
then order $W_B \leq n$, (5) for $w\in W_F$ such that $w\cap \text{int }F
\neq \eset$, then either $w\in W_B$ or int $F\supset \bar w$, and (6) if
$H\subset W_F$ and $N[H] = \ds{\bigcap_{h\in H}} h \neq \eset$, then $N[H]$
meets $F$.

\ni Proof.  For each $x\in F$, $\phi^{-1}\phi(x) = x$.  Since $\phi$ is an
open and closed mapping, it follows from [19; pp\. 78-80] that there is a
finite open irreducible covering $G_1$ of $F$ such that (1) for $u\in G_1$,
$\phi^{-1}\phi(u) = u$ (that is, $G_1$ is saturated) and (2) $G_1$ star
refines $W_2$.   

Now, dim $B = k \leq n-1$ and $B = \ds{\bigcup^{k+1}_{i=1}} K_i$ where dim
$K_i \leq 0$ [47].  There is a finite open refinement $\beta^i$ of $G_1$ which
covers $K_i$ where $\beta^i = \{v^i_1,v^i_2,\cdots,v^i_{n_i}\}$ such that
$v^i_s \cap v^i_t = \eset$ for $s\neq t$.  Let $G_2 = \bigcup \beta^i$.
Thus, order $G_2 \leq k + 1 \leq n$.  Let $G^*_2$ denote the union of the
elements of $G_2$.  Let $P$ denote an open set such that $P\supset B$ and
$G^*_2 \supset \bar P$.  Cover $F - G^*_2$ with a finite open refinement
$G_3$ of $G_1$ of order $\leq n+1$ such that if $g\in G_3$, then $\bar g
\cap \bar P = \eset$.  

Let $G_2 = \{g_1,g_2,\cdots,g_m\}$ and $G_3 = \{h_1,h_2,\cdots,h_s\}$.  By
[47; pp\. 11-33, 48; pp\. 133-135] 
there is an open refinement $G_4$ of $G_2 \cup G_3$ which covers
$F$ such that (1) $G_4 = \{u_1,u_2,\cdots,u_m;v_1,v_2,\cdots,v_s\}$, (2) for
each $i$, $u_i \subset g_i$ and $v_i\subset h_i$, and (3) order $G_4 \leq
n+1$.  Since $v_i \subset h_i$, $\bar v_i \cap \bar P = \eset$.  If $G'_2 =
\{u_1,u_2,\cdots,u_m\}$, then order $G'_2 \leq n$.  There is no loss of
generality in assuming that $G_4$ is irreducible and that if $H \subset G_4$ and
the nucleus, $N[H]$, is non empty, then $N[H]$ meets $F$.  To see the
latter, order the elements of $G_4$ as $U_1,U_2,\cdots,U_k$, $k = m+s$.  If
$E^r_j$ is a simplex of $G_4$, that is, in the nerve of $G_4$, $N(G_4)$,
then denote its nucleus by $N^r_j$.  If
$N^r_j$ meets $F$, then let $p^r_j \in F\cap N^r_j$; otherwise, let $p^r_j$
denote any point of $N^r_j$.  Let $N = \{p^r_j\mid p^r_j$ is chosen for
$N^r_j$ as above$\}$.  For each $i$, $N\cap U_i$ is a finite set in $U_i$.  By a
standard procedure, shrink the elements of $G_4$ to obtain an open covering $G'_4 =
\{U'_1,U'_2,\cdots,U'_{m+s}\}$ which covers $F$ such that for each $i$, $N\cap U_i
\subset U'_i \subset \bar U'_i \subset U_i$.   This can be accomplished as follows:
Replace each $U_i$ by an open set $U'_i$ such that (1) $U_i\supset \bar U'_i$ and (2)
$U'_i \supset (N\cap U_i) \cup \left(\bar P - \ds{\sum^{i-1}_{j=1}}U'_j -
\ds{\bigcup^k_{j=i+1}}U_i\right)$.  See [48:  pp\. 133-134].   

Suppose that a simplex ${E^r_i}' = U'_0 U'_1\cdots U'_r$ in the nerve of $G'_4$ has a
nucleus $Q^r_i$ that does not meet $F$.  Then $x\in \bar Q_i^r \cap F$
would imply $x\in U'_0 \cap U'_1 \cap \cdots\cap U'_r = N^r_i$ since $\bar
Q^r_i \subset \bar U'_0 \cap \cdots \cap \bar U'_r \subset N^r_i$.  But,
$F\cap N^r_i \neq \eset$, implying that $p^r_i \in U'_0 \cap \cdots 
U'_r = Q^r_i$, that is, $Q^r_i \cap F \neq \eset$. Consequently, if a
nucleus $Q^r_i$ fails to meet $F$, then $\bar Q^r_i \cap F = \eset$.
Finally, let $Q$ be an open set containing $F$ such that $Q$ meets no $\bar
Q^r_i$ that fails to meet $F$.  Let $G''_4 = \{U''_1,U''_2,\cdots,U''_{m+s}\}$
where $U''_i = U'_i \cap Q$ for $1\leq i\leq m+s$.  
If $H \subset G'_4$ and $N[H] \neq \eset$, then $N[H]
\cap F \neq \eset$ [cf\. 48, p\. 134].  The order of $G_4 \leq n+1$ implies
that order $G'_4 \leq n+1$ since $U_i \in G_4$ contains $U''_i \in G'_4$ for
each $i$, $1\leq i \leq m + s = k$.  Finally, let $W_F = G'''_4 =
\{U'''_1,U'''_2,\cdots,U'''_{m+s}\}$ where $U'''_i = \bigcup\{\phi^{-1}\phi(x) \mid
U^{\prime\prime}_i \supset \phi^{-1}\phi(x)\}$ is an open inverse set contained in
$U^{\prime\prime}_i$ ($U^{\prime\prime\prime}_i$ is open since $\phi$ is an open
and closed mapping and $\phi^{-1}\phi(U^{\prime\prime\prime}_i) =
U^{\prime\prime\prime}_i)$.  

Let $W_B = \{w\mid w \in W_F$ and
$w\cap B \neq \eset$ where $B = BdF\}$.  Since $W_B =
\{U^{\prime\prime\prime}_1,U^{\prime\prime\prime}_2,\cdots,U^{\prime
\prime\prime}_m\}$, $U^{\prime\prime\prime}_i\subset g_i \in G_2$, and order $G_2
\leq n$, then order $W_B \leq n$.  
It should be clear that $W_F$ has all properties claimed in
Lemma 3.  

Next, we extend $W_F$ to a {\it special covering} $V$ of $M$ (defined below) by
covering $M - \ds{\bigcup_{w\in W_F}}w$ in a special way.  We use the
following lemma.

\ni Lemma 4.  Suppose that $\phi \in L(M,p)$.  Then there exists a 
finite open covering $R$ of $Y=M/A_p$ such that (a)  if
$y\in\phi(F)$ where $F = F_\phi$,
then there is $r\in R$ such that $y\in r = \phi(w)$ for some
$w\in W_F$, (as described in Lemma 3),
(b) if $\phi(W_F) = \ds{\bigcup_{w\in W_F}}\phi(w)$ and $y\in Y
- \phi(W_F)$, then there is $r\in R$ such that $y\in r$, $r \in Q$ where $Q$ is the
basis in The Standing Hypothesis, 
$\bar r
\cap \phi(F) = \eset$, $\phi^{-1}(r) = r_1 \cup r_2 \cup \cdots \cup r_q$, $q
= p^t$ for some natural number $t$,
such that for each $i = 1,2,\hdots,q$, $r_i$ is a
component of $\phi^{-1}(r)$, $r_i$ maps onto $r$ under
$\phi$, $\overline{r_i} \cap \overline{r_j} = \eset$ for $i\neq j$, and
$r_i$ is homeomorphic to $r_j$ for each $i$ and $j$ with a homeomorphism compatible
with the projection $\phi$ (indeed, there is an
orientation preserving homeomorphism, an element of $A_p$, which takes $r_i$
onto $r_j$), (c)
$R$ is irreducible, (d) $V'=\{\bar c\mid c$ is a component of $\phi^{-1}(r)$
for some 
$r\in R\}$ such that $r \neq \phi(w)$ for any $w\in W_F\} \cup \{\bar w \mid w
\in W_F\}$ star refines $W_2$, and $V = \{v\mid \bar v \in V'\}$, and (e) if $r_x
\in R$ and $r_x\neq \phi(w)$ for any $w\in W_F$, $r_y \in R$ and $r_y \neq 
\phi(w)$ for any $w\in W_F$, $r_x\cap r_y \neq \eset$, $\phi^{-1}(r_x)$ consists
of exactly $p^{m_x}$ components, $\phi^{-1}(r_y)$ consists of exactly $p^{m_y}$
components, and $m_x \geq m_y$, then each component of $\phi^{-1}(r_y)$ meets
exactly $p^{m_x-m_y}$ components of $\phi^{-1}(r_x)$. 

\ni Proof.  Obtain $W_F$ using Lemma 3.  
Since $Y - \phi(W_F)$ is compact, use Lemma 1 to obtain 
a finite irreducible covering $R'$ of $Y - \phi(W_F)$ of sets $r$
satisfying the conditions of the lemma such that $R'$ star refines
$\{\phi(u) \mid u\in W_2\}$.  Property (e) of the conclusion of Lemma 4 is satisfied
by using the compactness of $Y$ and choosing $R'$ such that each $r\in R'$ has
sufficiently small diameter and $r\in Q$ (the basis in The Standing
Hypothesis). Let $R = R' \cup \{\phi(w)\mid w\in W_F\}$.
The lemma is established.

The irreducible finite open covering $V^1 = V$ of $M$ (which contains $W^1_F =
W_F$) generated by the irreducible finite open covering $R$ of $Y$ in Lemma 4 is
just the first step in establishing Lemma 5 below.  

\ni Lemma 5.  There are sequences $\{V^m\}$ and $\{U^m\}$ of finite open
coverings of $M$ cofinal in the collection of all open coverings of $M$ such that
(1) $V^{m+1}$ star refines $U^m$, (2) $V^1$ star refines $W_2$ of Lemma 2, (3)
$U^m$ star refines $V^m$, (4) order $U^m = n+1$, (5) $V^m$ is generated by a
finite open covering $R^m$ of $Y = M/A_p$, (6) $W^m_F$ is the subcollection of
$V^m$ which covers $F = F_\phi$, (7) $V^m$, $W^m_F$, and $R^m$ have the
properties stated in Lemma 4 where $R^m$ replaces $R$, $V^m$ replaces $V$, and
$W^m_F$ replaces $W_F$, (8) $\{$mesh $V^m\} \to 0$, (9) there are projections
$\pi_m : V^{m+1} \to V^m$ such that (a) $\pi_m = \beta_m\alpha_m$ where $\alpha_m
: V^{m+1} \to U^m$ and $\beta_m : U^m \to V^m$, (b) $\pi_m$ takes each
distinguished family $\big\{f^m_{ij}\big\}^{t^{m+1}_i}_{j=1}$ in $V^{m+1}$
(defined in a manner like those defined for $V^1$ and $V^2$ below) onto a
distinguished family $\big\{f^m_{sj}\big\}^{t^m_s}_{j=1}$ in $V^m$, 
and (c) $\pi_m$
extends to a simplicial mapping (also, $\pi_m$) of $N(V^{m+1})$ into
$N(V^m)$ such that if $\delta^n$ is an $n$-simplex in $N(V^{m+1})$ and if
$\pi_m(\delta^n) = \sigma^n$, an $n$-simplex in $N(V^m)$, then $N[\delta^n]
\cap F\neq \eset$ if and only if $N[\sigma^n] \cap F \neq \eset$.  (Also,
$\alpha_m$ and $\beta_m$ denote the extensions of $\alpha_m$ and $\beta_m$ to
simplicial mappings $\alpha_m : N(V^{m+1}) \to N(U^m)$ and $\beta_m : N(U^m) \to
N(V^m)$ where $\pi_m = \beta_m\alpha_m$.)  

NOTE: As stated in The Standing Hypothesis, $Y$
has a countable basis $Q = \big\{B_i\big\}^\infty_{i=1}$ such that (a) for
each $i$, $B_i$ is connected and uniformly locally connected and (b) if $H$
is any subcollection of $Q$ and $\ds{\bigcap_{h\in H}}h\neq \eset$, then
$\ds{\bigcap_{h\in H}}h$ is connected and uniformly locally connected (a consequence
of a theorem due to Bing and Floyd [50]).  Each $r\in R^m - \{\phi(w)\mid
w\in W_F\}$ can be chosen from
$Q$.

The proof of Lemma 5, although straightforward, is long and tedious.  The
existence of $V = V^1$ in Lemma 4 (which star refines $W_2$)
generated by $R = R^1$ is an initial step of a
proof using mathematical induction.  Additional first steps are described below.
These should help make it clear how the induction is completed to obtain a proof
of Lemma 5.  

Call the collection $\{r_1,r_2,\cdots,r_q\}$ consisting of all
components of $\phi^{-1}(r)$ a {\it distinguished}
family in $V$ (defined in Lemma 4)
generated by $r\in R - \phi(W^1_F)$ where $\phi(W^1_F) = \{\phi(w)\mid w\in W^1_F\}$
where $R$ satisfies Lemma 4.  
The finite open covering $V$ can be partitioned into either the
subcollections $\{r_1,r_2,\cdots,r_q\}$ consisting of the components of
$\phi^{-1}(r)$ for some $r\in R$ with $r\neq \phi(w)$ for any $w\in W_F^1$ and
$q = p^t > 1$ or collections of singletons $\{w\}$ where $w\in W_F^1$
which are defined to be {\it distinguished}
families generated by $r\in R$.  Thus, $V$ is generated by $R = R^1$.

Let $V^1 = V$.  
Clearly, $V^1$ star refines $W_2$.  
Observe that $V^1_F = \{v\mid v\in V^1$
and $v\cap F \neq \eset\} = W_F$ and that $V^1_B = \{v\mid v\in V^1_F$ and
$v\cap BdF \neq \eset\} = W_B$.  Also, order $V^1_B \leq n$ and order $V^1_F
\leq n+1$.  Note that order $V^1$ may be larger than $n+1$ since if $\phi$
is the orbit mapping of an effective action by a $p$-adic transformation
group, then dim $Y = n + 2$ or $\infty$ [22].

The covering $V^1=V$ is defined to be a {\it special covering of $M$ w.r.t}\.
$\phi$ generated by $R$. Of course, $\phi$ is fixed throughout
this discussion as in the statements of Lemmas 3 and 4.  
Observe that it follows from Lemma 1, 
that if $\big\{f^1_{kj}\big\}^q_{j=1}$ and $\big\{f^1_{mj}\big\}^s_{j=1}$ are
two non singleton families (those containing more than one element) in $V^1$
such that for some $i$ and $t$, $f^1_{ki} \cap
f^1_{mt} \neq \eset$, then for each $j$, the number of elements of
$\big\{f^1_{mj}\big\}^s_{j=1}$ which have a non empty intersection with
$f^1_{kj}$ is a constant $c_k$ and for each $j$, the number of elements of
$\big\{f^1_{kj}\big\}^q_{j=1}$ which have a non empty intersection with
$f^1_{mj}$ is a constant $c_m$ where $c_k = p^b$, $b\geq 0$, and $c_m = p^d$,
$d\geq 0$.  A singleton family (which contains exactly one element) either meets
each member of a non singleton family or meets no member of a non singleton
family.

\cl{Construction Of $U^1$ Of Order $n+1$ Which Refines $V^1$}

The reason that the sequence $\{U^m\}$ is constructed is to prove
(using the definitions of $\alpha_m$, $\beta_m$, and $\pi_m =
\beta_m\alpha_m$) that the inverse limit of the $n^{\text{th}}$
simplicial homology of the $n$-skeleta of the nerve of $V^m$ is $Z_p$
which permits the application of $\sigma$ (defined below) to actual
$n$-cycles.  The operator $\sigma$ can not be applied (as defined) to
elements of a homology class. 

{\it The next step is to describe a special refinement $U^1$
of $V^1$} which has order $n+1$ and other crucial properties.  
First, construct $U^1_1$.  List the {\it non degenerate} distinguished families of $V^1$
as $F^1_1,F^1_2,\cdots,F^1_{n_1}$ where $F^1_i = 
\big\{f^1_{ij}\big\}^{t^1_i}_{j=1}$ where $t^1_i = p^{c_i}$.
Recall that $f^1_{ij}$ is homeomorphic to $f^1_{it}$ for each $i$, $j$, and $t$
that makes sense.  Since $R^1$ (which generates $V^1$) is irreducible, it
follows that if $f^1_{ij}\in F^1_i$ and $f^1_{st} \in F^1_s$ where $F^1_i$
and $F^1_s$ are distinguished families in $V^1$ with $i \neq s$, then
$f^1_{ij}$ and $f^1_{st}$ are {\it independent}, that is, $f^1_{ij}
\not\supset f^1_{st}$ and $f^1_{st} \not\supset f^1_{ij}$.  This ordering may be
changed below.

For each $i$, $1\leq i \leq n_1$, choose a closed and connected
subset $K_i$ in $\phi(f^1_{i1}) = r_i \in R$ 
where $F^1_i = \big\{f^1_{ij}\big\}^{t^1_i}_{j=1}$ such that
(1) $K_{ij} = \phi^{-1}(K_i) \cap f^1_{ij}$, $K_{ij}$ is homeomorphic to
$K_{is}$ for any $s$ and $j$ that makes sense, (2) 
$C = \{\text{int }K_{ij}\mid 1\leq i \leq n_1$ and $1
\leq j \leq t^1_i\}$ covers $M - (W_F)^*$, (3) $K = \{\text{int }K_i\mid 1\leq i\leq
n_1\}$ covers $Y - \ds{\bigcup_{w\in W_F}}\phi(w)$, and (4) $K_{ij}$ is
connected, $1\leq j\leq t^1_i$.
To see that this is possible, choose a closed subset $A_i$ of $r_i$, $1\leq i\leq n_1$,
such that $A = \{A_i\mid 1\leq i\leq n_1\}$ covers $Y - \ds{\bigcup_{w\in W^1_F}}\phi(w) =
Y'$ [47; 49].  Note that there exists
a natural number $k$ such that $\{A_i = r_i - N_{\frac{1}{k}}(\partial r_i)\mid 1\leq
i\leq n_1\}$ covers $Y'$.  To see this, suppose that for each $k$, there is $x_k \not\in
\ds{\bigcup^{n_1}_{i=1}}(r_i - N_{\frac{1}{k}}(\partial r_i))$.  Since $Y'$ is compact,
there is a subsequence $\{x_{n(k)}\}$ of $\{x_k\}$ which converges to $x\in r_q$ for some
$q$.  There is some $m$ such that $x\in r_q - N_{\frac{1}{m}}(\partial r_q)$ and $x\in$
interior$(r_q - N_{\frac{1}{m+1}}(\partial r_q))$ which leads to a contradiction.
Choose $p_i \in r_i$ and $m$ sufficiently large such that $K_i$ is the component of
$r_i - N_{\frac{1}{m}}(\partial r_i)$ which contains $p_i$.  For each $m$, let
$C_m(p_i)$ be the component of $r_i - N_{\frac{1}{m}}(\partial r_i)$ which contains
a fixed $p_i \in r_i$ (with $m$ large enough).  It will be shown that
$\ds{\bigcup^\infty_{m=1}}C_m(p_i) = r_i$.  Suppose that there is $q\in r_i$ such
that $q\not\in \ds{\bigcup^\infty_{m=1}}C_m(p_i)$.  Since $r_i$ is uniformly
locally connected and locally compact ($\bar r_i$ is a Peano continuum), there is a
simple arc $p_iq$ from $p_i$ to $q$ in $r_i$.  Consequently, for some $m$, $r_i -
N_{\frac{1}{m}}(\partial r_i) \supset p_iq$ which is in $C_m(p_i)$.  This is
contrary to the assumption above.  Hence, $\ds{\bigcup^\infty_{m=1}}C_m(p_i) =
r_i$.  For $p_i \in r_i$ fixed as above, choose $p_{ij} \in \phi^{-1}(p_i) \cap
f^1_{ij}$ for $1\leq j\leq t^1_i$.  Let $C^i_m(p_{ij})$ be the component of
$\phi^{-1}(C_m(p_i))$ which is in $f^1_{ij}$ and contains $p_j$.  It will be shown
that $\ds{\bigcup^\infty_{m=1}}C^i_m(p_{ij}) = f^1_{ij}$.  If this is false, then
there is a $q_j \in f^1_{ij} - \ds{\bigcup^\infty_{m=1}}C^i_m(p_{ij})$.  Since
$f^1_{ij}$ is ulc and locally compact, there is a simple arc $p_{ij}q_j$ from $p_{ij}$ to
$q_j$ in $f^1_{ij}$.  Now, $r_i \supset \phi(p_{ij}q_j)$.  For $m$ large enough,
$C_m(p_i) \supset \phi(p_{ij}q_j)$ and some component of $\phi^{-1}\phi(p_{ij}q_j)$
contains $p_{ij}q_j$ and lies in $C^i_m(p_{ij})$. This is contrary to the assumption
above.  Hence, $\ds{\bigcup^\infty_{m=1}}C^i_m(p_{ij}) = f^1_{ij}$.  Choose $m$
sufficiently large that condition (2) above is satisfied where $K_i =
C_m(p_i)$ and (4) int $C_{m}(p_i)\supset A_i$.  By [19], $\phi^{-1}(C_m(p_i))$ consists
of a finite number of components each of which maps onto $K_i$ under $\phi$.  Choose $m$
sufficiently large that the components of $\phi^{-1}(K_i)$ are
$\big\{K_{ij}\big\}^{t^1_i}_{j=1}$ and conditions (1) - (4) are satisfied.  

Next, shrink the elements of $W^1_F=W_F$ as follows:  Order $W^1_F$ as
$w^1_1,w^1_2,\cdots,w^1_{f_1}$ and choose a natural number $k$ such that $\{K_{n_1+i}
= \phi(w^1_i) - N_{\frac{1}{k}}(\partial \phi(w^1_i))\mid 1\leq i\leq f_1\}$ has the
property that $\{\text{int }K_{n_1+i}\mid 1\leq i\leq f_1\}$ covers $\ds{\bigcup_{w\in
W^1_F}}\phi(w) - \ds{\bigcup^{n_1}_{i=1}}\text{int }K_i$.  
 
By a Theorem (Nagami and Roberts) and its Corollary [49; pp\. 90-91], a metric
space $X$ has dim $X\leq n$ if and only if $X$ has a sequence $\{G_i\}$ of open
coverings of $X$ such that (1) $G_{i+1}$ refines $G_i$ for each $i$, (2) order $G_i
\leq n+1$ for each $i$, and (3) mesh $G_i < \frac{1}{i}$.
If $X$ is a manifold, then the elements of $G_i$ can be chosen to be connected and
uniformly locally connected.  If $X$ is a triangulable $n$-manifold, then it is easy to see
this by using barycentric subdivisions of a triangulation of $X$.  Choose such a sequence
$\{G_i\}$ of open coverings of $M$ such that $G_1$ star refines 
$K = \{\text{int }K_i\mid 1\leq i\leq n_1\}$.

The next step is to show how to choose some $G_i$ from which $U^1_1$ is chosen and later
modified to give $U^1$ with the desired properties.

Let $\epsilon$ be the Lebesque number of the covering $K$.
Choose $t$ such that (1) mesh $G_t < \frac{1}{t}$, (2) if $g\in G_t$, then
diam $\phi(g) < \frac{\epsilon}{8}$, and (3) if $y\in Y$ and $G(y) = \{g\mid g\in G_t$
and $y\in \phi(g)\}$,
then there is $s$, $1\leq s\leq n_1+n_2$, such that int $K_s \supset
\ds{\bigcup_{g\in G(y)}}\overline{\phi(g)}$.
(Note that (3) follows from the
choice of $\epsilon$, $t$ and the fact that $G_t$ star refines $K$.).

Choose $U^1_1 \subset G_t$ such that $U^1_1$ is an
irreducible finite covering of $M$.  
Since $G_t$ is an open covering of $M$ with connected open sets, mesh $G_t
< d$, and order $G_t \leq n+1$, $U^1_1$ has the following properties:
\spitem{(1)}  order $U^1_1 \leq n+1$ and $U^1_1$ star refines $V^1$, 
\spitem{(2)}  if $u\in U^1$, then $u$ is connected, and
\spitem{(3)}  if $y\in Y$ and $G(y) = \{g\mid g\in U^1_1$ and $y \in
\phi(g)\}$,
then there is $s$, $1\leq s \leq n_1+n_2$, such that int $K_s
\supset \ds{\bigcup_{g\in G(y)}}{\overline{\phi(g)}}$.
\vskip 5pt

\cl{Construction Of $U^1$ Which Refines $U^1_1$ In A Special Way}

For each $y\in Y$, let $Q(y) =\{\phi(u)\mid u \in U^1_1$
and $y\in \phi(u)\}$.  There are at most a finite number of such sets
distinct
from each other.  Order these sets as $Q_1,Q_2,\cdots,Q_{m_1}$ such that for $i < j$, $Q_i
\neq Q_j$ and card $Q_i \geq$ card $Q_j$.  Let $O_i = \cap Q_i - \ds{\bigcup_{j <
i}}(\overline{\cap Q_j})$
where $\cap Q_i = \{x\mid x \in \phi(u)$ for each $\phi(u) \in Q_i$,
$1\leq i\leq m_1\}$.  For each $i$, $1\leq i\leq m_1$, let $B_{i\phi(u)} = (\partial
O_i)\cap \partial \phi(u)$ where $\phi(u)\in Q_i$ and $(\partial O_i) \cap \partial
\phi(u) \neq \eset$.
There are at most a finite number of such non empty sets distinct
from each other.  Let $B_1,B_2,\cdots,B_{m_2}$ denote all those sets distinct from each
other.  Let $B = \ds{\bigcup^{m_2}_{i=1}}B_i$.  For each $y\in B$, let $D(y) = \{B_t \mid
y\in B_t\}$.  There are at most a finite number of such non empty sets distinct from each
other.  Order these as $D_1,D_2,\cdots,D_{m_3}$ such that if $i < j$, then $D_i \neq D_j$
and card $D_i \geq \text{card }D_j$.  For each $i$, $1\leq i\leq m_2$, there is $u\in
U^1_1$ and a closed subset $C_{iu}$ of $\partial u$ such that $\phi(C_{iu}) = B_i$.

There is a finite open irreducible saturated covering $W^2_F$ which consists of
open inverse sets, $W^2_F$ covers $F = F_\phi$, star refines
both $W_F = W^1_F$ in Lemma 3 and $C$, and $W^2_F$ has properties (2)-(6)
attributed to $W_F$ in Lemma 3, and for each $i$, $1\leq i\leq n_1$, and $w\in W^2_F$,
$\phi^{-1}(r_i) \cap \bar w = \eset$ where $r_i \in R^1 - \phi(W^1_F)$.  

It follows from the definition that 
$\phi^{-1}(B)$ is closed and contains no open set.  Hence, dimension $\phi^{-1}(B) 
\leq n-1$.  For each $y\in Y$, choose $0 < \epsilon_y <
(\frac{1}{8})\min\{\rho(y,B_i)\mid y\not\in B_i$, $1\leq i\leq m_2\} \cup \rho(\cap
D_i,\cap D_j)\mid(\cap D_i)\cap (D_j) = 0\}$ and cover $B$ with a finite irreducible
collection $E$ of $\epsilon_y$-neighborhoods for $y\in B$.  Choose $t$ sufficiently
large such that if $g\in G_t$ and $\phi(g)\cap B \neq \eset$, then there is $e\in E$
such that $e\supset \overline{\phi(g)}$.  Let $H'$ be a finite irreducible collection
of open sets which refines $G_t$, covers $\phi^{-1}(B) - (W^2_F)^*$, refines $U^1_1$,
and order $H' \leq n$ (See [48; pp.\ 133-347]).  In addition choose $H'$ such that 
$H'$ can be extended to a finite
irreducible cover $G$ of $M$ (i.e., $G\supset H')$ 
such that order $G\leq n+1$, $G$ refines
$U^1_1$, and if $g\in G-H'$, then $\bar g \cap \phi^{-1}(B) = \eset$.  Now, cover
$\phi^{-1}(B)$ with a finite irreducible collection $H$ of open sets such that (a) $H$
refines $H'$, (b) order $H\leq n$, (c) order $(H\cup H') \leq n+1$ (if $N[\delta^{n-1}]$ is
the nucleous of an $(n-1)$-simplex in the nerve, $N(H')$, of $H'$, then one and only one element
of $H$ contains $N[\delta^{n-1}] \cap \phi^{-1}(B))$, and (d) if $O = \{g\cap \phi^{-1}(O_i)
\mid g\in G$ and $1\leq i\leq m_1\}$ (recall that $\phi^{-1}(B)\cap\phi^{-1}(O_i) = \eset$),
then $U^1_2 = H\cup O$ has the following properties:
(1) $U^1_2$ star refines $U^1_1$, (2)
order $U^1_2=n+1$, (3) if $(\cap D_i) \cap (\cap D_j) = \eset$, $u\in U^1_2$, $v\in U^1$,
$\phi(u) \cap (\cap D_i) \neq \eset$, and $\phi(v) \cap (\cap D_j) \neq \eset$, then
$\overline{\phi(u)} \cap \overline{\phi(v)} = \eset$, (4) if $B_i \cap B_j = \eset$,
$u\in U^1_2$, $v\in U^1_2$, $\phi(u) \cap B_i \neq \eset$, and $\phi(v) \cap B_j \neq \eset$,
then $\overline{\phi(u)} \cap \overline{\phi(v)} = \eset$ (a consequence of (3)), 
(5) if $u\in U^1_2$, $\phi(u) \cap (\cap D_i) \neq \eset$, and $\phi(u) \cap (\cap D_j)
\neq \eset$, then $(\cap D_i) \cap (\cap D_j) \neq \eset$, (6) if $u\in U^1_2$ and
$\phi(u) \cap B= \eset$ where $B = \ds{\bigcup^{m_2}_{j=1}}B_j$, then $O_i \supset
\phi(u)$ for some $i$, $1\leq i\leq n_1$ (it follows that $\phi(h) \supset
\overline{\phi(u)}$ for all $h\in U^1_1$ such that $\phi(h) \cap \phi(u) \neq \eset)$,
and (7) if $u\in U^1_2$ and $\phi(u) \cap (\cap D_t) \neq \eset$ for the smallest $t$, then
$\phi(h) \supset \phi(u)$ for all $h\in U^1_1$ such that there is no $B_x \in
D_t$ with the property that $\partial\phi(h) \supset B_x$ and $\phi(h) \cap \phi(u) \neq
\eset$. 
It is not difficult to see that $U^1_2$ has properties (1)-(7), 
$U^1 = U^1_2 \cup W^2_F$ has order $n+1$, and
$w\in W^2_F$, then $u\not\supset
w$ for any $u\in U^1_2$.  

Observe that if $u\in U^1_2$ and $\phi(u) \cap B = \eset$, then for some $i$, $O_i \supset
\phi(u)$ since there is some $e\in E$ such that $e\supset \overline{\phi(u)}$.

\cl{Construction of $V^2$ Which Refines $U^1$}

Next, construct $V^2$.  
For each $y\in Y$, let $G(y) =
\{u\mid u\in U^1_1$ and $y\in \phi(u)\}$.   There is some $s$, $1\leq
s\leq n_1+n_2$, such that $K_s \supset \ds{\bigcup_{u\in G(y)}}{\overline{\phi(u)}}$.

For $y\in B - \phi(W^2_F)^*$, 
choose $r_y \in Q$ (the basis for $Y$ described above)
such that (1) $\ds{\bigcap_{u\in G(y)}}\phi(u) \supset \bar r_y$, (2) if $U(y) = \{u\mid
u\in U^1$ (not $U^1_1$) and $y\in \phi(u)\}$, then $\ds{\bigcap_{u\in U(y)}}\phi(u)
\supset \bar r_y$, (3)  diam $r_y <
(\frac{1}{8})\min\{ \rho(y,\partial\phi(v))\mid v\in U^1_1$ and $y\not\in
\partial\phi(v)\}$, (4) $y\in r_y$, (5) if $w\in W^2_F$, then $\phi(w) \not\supset r_y$,
(6) $\phi^{-1}(r_y) = r_{y1} \cup r_{y2} \cup \cdots \cup r_{yq}$, $q = p^{t_q}$
where $t_q \geq 1$, $r_{yi}$ maps onto $r_y$ under $\phi$, $\bar r_{yi} \cap \bar
r_{yj} = \eset$ for $i\neq j$, and $r_{yi}$ is homeomorphic to $r_{yj}$ for each $i$ and
$j$ with a homeomorphism compatible with the projection $\phi$ (indeed, there is an
element of $A_p$ which takes $r_{yi}$ onto $r_{yj}$),
and (7) for each $i$, there is $u\in
U^1$ such that $u \supset \bar r_{y_i}$.  See the proof of Lemma 4 for the construction.  

Let $R^2_1$ denote a finite irreducible collection of such sets $r_y$ which covers
$B - \phi(W^2_F)^*$.   If $y\in Y$ and
$y\not\in (R^2_1)^* - \phi(W^2_F)^*$,  
then choose $r_y$ satisfying (1) - (7) above such that
$\bar r_y \cap B = \eset$ and let
$R^2_2$ denote a finite irreducible cover of $(R^2_1)^* - \phi(W^2_F)^*$
consisting of such $r_y$.  Let $R^2 = R^2_1 \cup R^2_2 \cup
\phi(W^2_F)$ which
is an irreducible cover of $Y$.  Let $V^2 =
\{c\mid c$ is a component of $\phi^{-1}(r_{y_i})$ for some $i$, where $r_{y_i} \in R^2_1
\cup R^2_2\}\cup W^2_F$ which  
is an irreducible cover of $M$ that star refines $V^1$ and $U^1$.
The collection of components, $\big\{f^2_{ij}\big\}^{t^2_i}_{j=1}$, of
$\phi^{-1}(r_{y_i})$, $r_{y_i} \in R^2_1 \cup R^2_2$ 
is a non degenerate distinguished family in $V^2$ whereas each $w\in
W^2_F$ where $w = \phi^{-1}\phi(w)$ is a singleton distinguished family in $V^2$.
 
\cl{Definitions Of $U^1$, $\alpha_1$, $\beta_1$, And $\pi_1 = \beta_1\alpha_1$}

Case (1): $y_i \in B  -
\phi(W^2_F)^*$ and $y_i \in r_{y_i} \in R^2_1$.
Let $F^2_i = \big\{f^2_{ij}\big\}^{t^2_i}_{j=1}$ be any non degenerate distinguished
family in $V^2$ generated by $r_{y_i}$ in $R^2_1$
in this case.  Let $c_i = \min\{t\mid y_i \in \cap D_t\}$.  Let $s_i =
\min\{s\mid \text{int }K_s\supset {\overline{\phi(u)}}$ for all $u\in U^1_1$ such that $u \in
G(y_i)\}$. Choose $F^1_{s_i} =
\big\{f^1_{s_ij}\big\}^{t^1_{s_i}}_{j=1}$ for $F^2_i$. For each $j$, $1\leq j\leq t^2_i$,
choose $U_{ij} \in H\subset U^1$ such that $U_{ij} \supset f^2_{ij}$. There is a unique
$z_{ij}$, $1\leq z_{ij} \leq t^1_{s_i}$, such that $f^1_{s_iz_{ij}} \supset U_{ij}
\supset f^2_{ij}$. Let $\alpha_1(f^2_{ij}) = U_{ij}$, $\beta_1(U_{ij}) =
f^1_{s_iz_{ij}}$, and $\pi_1 =\beta_1\alpha_1$.

Case (2):  $y_i \not\in B$, $y_i \not\in \phi(W^2_F)^*$,
and $r_{y_i} \in R^2_2$, in this case, generates $F^2_i =
\big\{f^2_{ij}\big\}^{t^2_i}_{j=1}$ in $V^2$.  Let $s_i = \min\{s\mid \text{int }K_s \supset
\ds{\bigcup_{u\in G(y_i)}}{\overline{\phi(u)}}\}$.
Choose $F^1_{s_i} =
\big\{f^1_{s_ij}\big\}^{t^1_{s_i}}_{j=1}$ for $F^2_i$.  For each $j$, $1\leq j\leq
t^2_i$, choose $U_{ij} \in O\subset U^1_2$ 
such that $U_{ij} \supset f^2_{ij}$. There is a unique $z_{ij}$, $1\leq z_{ij} \leq
t^1_{s_i}$, such that $f^1_{s_iz_{ij}} \supset U_{ij} \supset f^2_{ij}$.  Let
$\alpha_1(f^2_{ij}) = U_{ij}$, $\beta_1(U_{ij}) = f^1_{s_iz_{ij}}$, and $\pi_1 =
\beta_1\alpha_1$.

Case (3): $y_i \in w\in W^2_F$.   
If $w\in W^2_F$, then $w\in U^1$.
Furthermore, there is $s_w = \min\{s\mid K_s\supset w,n_1\leq s\leq n_2\}$.  Choose
$w^1_{s_w}$ in $W^1_F$ such that $K_s$ is the shrinking of $w^1_{s_w}$.  Let $\alpha_1(w)
= w$ and $\beta_1(w) = w^1_{s_w}$.  

It will be shown that the mappings $\alpha_1$ and $\beta_1$ are well defined.

Case (A):  $w\in W^2_F$.  It should be clear that $\beta_1$ is well defined since the
choice of $u \in U^1$ such that $u\supset w$ is $u=w$.

Case (B):  Suppose that $\beta_1$ is not well defined and there exists $F^2_i$ and
$F^2_k$, two different non degenerate distinguished families in $V^2$ such that (a)
$s_i\neq s_k$ (if $s_i = s_k$, then $\beta_1$ is well defined), $F^1_{s_i}$ is chosen for
$F^2_i$, $F^1_{s_k}$ is chosen for $F^2_k$, and (b) $U_{ij} = U_{kt} \supset f^2_{ij} 
\cup f^2_{kt}$ where $F^2_i =
\big\{f^2_{ij}\big\}^{t^2_i}_{j=1}$, $F^2_k = \big\{f^2_{kj}\big\}^{t^2_k}_{j=1}$, for
some $j$, $1\leq j\leq t^2_i$, $U_{ij} \in U^1$ is chosen such that $U_{ij} \supset
\overline{f^2_{ij}}$, and for some $t$, $1\leq t\leq t^2_k$, $U_{kt} = U_{ij} \in U^1$ is
chosen such that $U_{kt}\supset \overline{f^2_{kt}}$ as described above.  
  
Case B(1):  $U_{ij} = U_{kt} \in O\subset U^1_2 \subset U^1$.  In this case, $y_i \not\in B$ and $y_k
\not\in B$.  Indeed, $y_i \in O_x$ and $y_k \in O_x$ for some $x$, $1\leq x\leq n_1$.
In this case, it follows from Property (6) of the properties of $U^1$ that for each
$u\in G(y_i)$, $y_k \in \phi(u)$, and for each $v\in G(y_k)$, $y_i \in \phi(v)$ since
$O_x \supset \phi(U_{ij}) =\phi(U_{kt})$ and $\cap Q_x\supset O_x$ as defined above.
Thus, $G(y_i) = G(y_k)$ and $s_i = s_k$ contrary to the assumption above.

Case B(2):  $U_{ij} = U_{kt} \in H\subset U^1$, $y_i \in B$, and $y_k
\in B$.  Recall that $c_i = \min\{t\mid y_i \in \cap D_t\}$
and $s_i = \min\{s\mid \text{int }K_s \supset
{\overline{\phi(u)}}$ for all $u\in U^1_1$ such that $u\in G(y_i)\}$.
Also, $c_k = \min\{t\mid y_k \in \cap D_t\}$ and $s_k =
\min\{s\mid \text{int }K_s\supset \overline{\phi(u)}$ for all $u\in U^1_1$ such that $u \in
G(y_k)\}$.
Furthermore, $U_{kt} \in H$ where $\phi(U_{kt}) \cap (\cap
D_{c_k}) \neq \eset$ and $U_{kt} = U_{ij} \in H$ such that
$\phi(U_{kt}) \cap (\cap D_{c_i}) \neq \eset$.  Now, $y_i \in \phi(U_{kt}) =
\phi(U_{ij})$ and $y_k \in \phi(U_{ij})$.
Now, $y_i \in \phi(u)$ for each $u\in G(y_i)$. Since $y_i \in \phi(U_{ij})$, it
follows by construction of $U^1$ that $\phi(h) \supset \phi(U_{ij})$ for all
$h\in U^1_1$ such that there is no $B_x \in D_{c_i}$ such that $\partial\phi(h) \supset
B_x$ and $\phi(h) \cap \phi(U_{ij}) \neq \eset$.  If there is $B_x \in D_{c_i}$ such that
$\partial\phi(u) \supset B_x$, then $y_i \in B_x$ and $y_i \not\in \phi(u)$.
Consequently, $\phi(u) \supset \phi(U_{ij}) = \phi(U_{kt})$ and
$y_k \in \phi(u)$.  Clearly, $y_k \in \phi(v)$ for each $v\in G(y_k)$.  It follows in
a similar way that $y_i \in \phi(v)$.  Thus, $G(y_i) = G(y_k)$ and
$s_i = s_k$ contrary to the assumption above.

Case (C):  Suppose that $\beta_1$ is not well defined and there exists $F^2_i$, a non
degenerate distinguished family in $V^2$ and $w\in W^2_F$, a singleton distinguished
family, such that $s_i \neq s_w$, $F^1_{s_i}$ is chosen for $F^2_i$, $w^1_{s_w} \in
W^1_F$ is chosen for $w$, and for some $j$, $1\leq j\leq t^2_i$, $U_{ij} = w\supset
f^2_{ij} \cup w$.  This contradicts the choice of $r_{y_i}$ such that $\phi(W^2_F)^* 
\not\supset r_{y_i}$.  Here, $\phi(U_{ij}) = \phi(w) \supset r_{y_i}$.  

It should be clear that $\alpha_1$ and $\beta_1$ are well defined.  

Clearly, $\beta_1$ is defined on $U^1$ since $U^1$ is irreducible and $V^2$ 
refines $U^1$.  Observe that $\pi_1$ maps
distinguished families onto distinguished families.

The first steps in the proof of Lemma 5 are complete.

With $V^i$ defined (as indicated for $i=1$ and 2), define $U^i$ in the
manner that $U^1$ is defined for $V^1$.  Use Lemma 4 to obtain a finite open
covering $R^{i+1}$ of $Y$ satisfying the conditions of Lemma 4 where $V^i$
replaces $W_2$, and $W^{i+1}_F$ replaces $W_F$,
and $R^{i+1}$ replaces $R$ such that $R^{i+1}$ generates a
special covering $V^{i+1}$ of $M$ having the properties similar to those
described above for $V^2$ w.r.t\. $U^1$ and $V^1$ but w.r.t\. $U^i$ and
$V^i$.  Similarly, define $\alpha_i$, $\beta_i$, and $\pi_i$ in the manner
that $\alpha_1$, $\beta_1$, and $\pi_1$ are described above.  
Extend $\alpha_i$, $\beta_i$, and $\pi_i$ in
the usual manner to the nerves $N(V^{i+1})$, $N(U^i)$, and $N(V^i)$,
respectively.

It should be clear that the proof of Lemma 5 can be completed using mathematical
induction and the methods employed above.

\cl{Orientation Of The Simplices In $N(V^m)$}

{\it Next, orient the simplices in $N(V^m)$ for each special covering}
$V^m$.  Recall that for each special covering $V^m$, there is associated a
covering $R^m$ of $Y$ which generates $V^m$.  Let $R = R^m$ and $V = V^m$.  
Suppose that $(v_0,v_1,v_2,\cdots,v_k) = \sigma^k$ is a
$k$-simplex in $N(R)$.  For each $i$, $0 \leq i \leq k$, $\phi^{-1}(v_i) =
v_{i1} \cup v_{i2} \cup \cdots \cup v_{it_i}$ where
$\{v_{i1},v_{i2},\cdots,v_{it_i}\}$ is the distinguished family determined
by $v_i$.  Consequently, $\sigma^k$ determines a distinguished family of
$k$-simplices in $N(\phi^{-1}(R))$ where $\phi^{-1}(R) = \{c \mid c$ is a
component of $\phi^{-1}(r)$ for $r\in R\}$.   
If $q\in N[\sigma^k]$, the nucleus or carrier of $\sigma^k$, then $\phi^{-1}(q)
\cap v_{ij} \neq \eset$ for each $i = 1,2,\cdots,k$ and each $j =
1,2,\cdots,t_i$.  The orientation of a $k$-simplex $\delta^k =
(v_{0j_0},v_{1j_1},\cdots,v_{kj_k})$ is to be that of $\sigma^k$ as
indicated by the given order of the vertices
$(v_{0j_0},v_{1j_1},\cdots,v_{kj_k})$ of $\delta^k$.  Since $\phi^{-1}(v_i)$
has $t_i$ components, there will be at least $t_i$ $k$-simplices in
$N(\phi^{-1}(R))$ which are mapped to $\sigma^k$ by the simplicial mapping 
$\phi^*
: N(\phi^{-1}(R)) \to N(R)$ induced by $\phi$.  This collection of
$k$-simplices is the {\it distinguished family determined by} $\sigma^k$
(more precisely, determined by $N[\sigma^k]$).  If $N[\sigma^k] \cap \phi(F)
= \eset$, then (1) $v_i \cap \phi(F) = \eset$ for some $i$, (2) the
distinguished family of $k$ simplices in $N(V^m)$ has cardinality $p^c$ for
some natural number $c$, and (3) $N[\sigma^k]$ is connected and ulc.  If
$N[\sigma^n]\cap \phi(F) \neq \eset$ where $\sigma^n$ is an $n$-simplex in
$N(R)$, then (1) int $\phi(F) \supset N[\sigma^n]$ and (2) $\sigma^n$
determines a distinguished family of one $n$-simplex $\delta^n$ in $N(V^m)$ and
int $F\supset N[\delta^n]$ by construction.    

\cl{Distinguished Families Of $n$-Simplices In $N(V^m)$}

The distinguished families $F^1_i = \big\{f^1_{ij}\big\}^{t^1_i}_{j=1}$ of members of
the covering $V^1$ generate distinguished families of $n$-simplices. That is, for
distinguished families $F^1_{k_i}$, $0\leq i \leq n$, in $V^1$ such that
$\ds{\bigcap^n_{i=1}}(F^1_{k_i})^* \neq \eset$ where $(F^1_{k_i})^* =
\ds{\bigcup^{t^1_i}_{j=1}}f^1_{k_ij}$, the $F^1_{k_i}$, $0\leq i\leq n$, generate
a distinguished family of $n$-simplices consisting of all $n$-simplices
$\{f^1_{k_0j_0},f^1_{k_1j_1},\cdots,f^1_{k_nj_n}\}$ such that (1) $f^1_{k_ij_i}
\in F^1_{k_i}$, $0\leq i\leq n$, and (2) $\ds{\bigcap^n_{i=0}}f^1_{k_ij_i} \neq
\eset$.  The number of $n$-simplices in this family is $\max\{t^1_{k_i} \mid 0
\leq i\leq n\}$.  As pointed out above, a distinguished family of $n$-simplices
is determined by an $n$-simplex $\sigma^n$ in $N(R)$ and such a family is a
lifting of $\sigma^n$ to $N(V^m)$.

Consider a distinguished family $S^1_k$ of $n$-simplices in
$N(V^1)$ as defined above such that a distinguished family $S^2_q$ of $n$-simplices
defined similarly in $N(V^2)$ using distinguished families in $V^2$ is mapped
onto $S^1_k$ by $\pi^*_1: N(V^2) \to N(V^1)$.  By construction, each $n$-simplex
in $S^1_k$ is the image of exactly $p^c$ $n$-simplices for fixed $c$, a
non negative integer, where card $S^2_q = p^{c_q}$, card $S^1_k = p^{c_k}$,
and $c = c_k - c_q$.  Of course, as a chain in $N(V^1)$, this is a trivial
$n$-chain using coefficients in $Z_p$.  

For each natural number $m$, define distinguished families $S^m_i$ of
$n$-simplices in $N(V^m)$ as described above.  Each such family
$S^m_i$ is the lifting of an $n$-simplex $\gamma^n_i$ in $N(R^m)$, that is,
$\phi : M \to M/A_p$ induces a mapping $\phi^* : N(V^m) \to N(R^m)$.  If
$\gamma^n_i$ is an $n$-simplex in $N(R^m)$, then $(\phi^*)^{-1}(\gamma^n_i)$
is the union of a distinguished family of $n$-simplices $S^m_i$ in $N(V^m)$.
The family $S^m_i$ is the lifting of $\gamma^n_i$ in $N(V^m)$.  Each member
of $S^m_i$ has the same orientation as $\gamma^n_i$.  

Observe that if $S^m_i$ and $S^m_j$ are two non degenerate distinguished
families of $n$-simplices such that $S^m_i
=\big\{\delta^n_{it}\big\}^{p^{c_i}}_{t=1}$, $S^m_j =
\big\{\delta^n_{jt}\big\}^{p^{c_j}}_{t=1}$, and $\delta^n_{it}$ shares and
$(n-1)$-face with $\delta^n_{js}$ for some $j$ and $s$, then for each $t$,
$1\leq t \leq p^{c_i}$, there is some $s$, $1\leq s \leq p^{c_j}$, such that
$\delta^n_{it}$ shares an $(n-1)$-face with $\delta^n_{js}$ and, conversely,
for each $s$, $1\leq s \leq p^{c_j}$, there is some $t$, $1\leq t \leq
p^{c_i}$, such that $\delta^n_{js}$ and $\delta^n_{it}$ share an
$(n-1)$-face.  If $c_i < c_j$, then for each $t$, $1\leq t\leq p^{c_i}$,
$\delta^n_{it}$ shares an $(n-1)$-face with exactly $p^{c_j-c_i}$ members of
$\delta^m_j$.  

If an $n$-simplex $\delta^n$ in $N(V^m)$ is such that $F_\phi \supset
N[\delta^n]$, the nucleus of $\delta^n$, then $\delta^n$ constitutes a
singleton distinguished family in $N(V^m)$ as described above.

\cl{The $n$-Skeleton, $N(V^m_n)$, Of $N(V^m)$, And Inverse Limit $H_n(V^m_n) \cong
Z_p$}
\sk1

Let $N(V^m_n)$ denote the $n$-skeleton of $N(V^m)$.  Recall that $V^m$
denotes the special $m^{\text{th}}$
covering. The special projections $\pi^*_m : N(V^{m+1}) \to N(V^m)$,
factors by $\alpha_m^* : N(V^{m+1}) \to N(U^m)$ and $\beta^*_m : N(U^m) \to
N(V^m)$ where $\pi^*_m = \beta^*_m\alpha^*_m$.

Let $H_n(V^m_n)$ denote the $n^{\text{th}}$ simplicial homology of
$N(V^m_n)$, the $n$-skeleton of $N(V^m)$.  
The coefficient group is always $Z_p$.
The following two lemmas give facts concerning the homology which will be
needed to finish the proof of the Theorem.

First, consider the following commutative diagram:
$$\matrix
\gets &H_n(V^m) &\overset {\beta^*_m} \to \longleftarrow &H_n(U^m)
&\overset {\alpha^*_m} \to \longleftarrow &H_n(V^{m+1}) &\overset
{\beta^*_{m+1}}\to \longleftarrow &H_n(U^{m+1}) &\gets & \\
& & & & & & & & & \\
 &{\ssize{\nu_m}} \uparrow &{\ssize{\beta^*_m}}\swarrow & &\nwarrow
{\ssize{\alpha^*_m}} &{\ssize{\nu_{m+1}}} \uparrow &{\ssize{\beta^*_{m+1}}}
\swarrow & &\nwarrow {\ssize{\alpha^*_{m+1}}} & \\
& & & & & & & & & \\
\gets &H_n(V^m_n) & &\overset {\pi^*_m} \to \longleftarrow & &H_n(V^{m+1}_n)
& &\overset{\pi^*_{m+1}} \to \longleftarrow &H_n(V^{m+2}_n)
&\gets\endmatrix$$

Here $\nu_m$ is the natural map of a cycle in $H_n(V^m_n)$ into its homology
class in $H_n(V^m)$.  The other maps are those induced by the projections
$\alpha_m$, $\beta_m$ and $\pi_m$.  The upper sequence, of course, yields
the {\v C}ech homology group ${\check H}_n(M)$ as its inverse limit.
Furthermore, it can be easily shown, using the diagram, that ${\check
H}_n(M)$ is isomorphic to the inverse limit $G = \varprojlim H_n(V^m_n)$, of
the lower sequence.  Specifically, $\gamma : {\check H}_n(M) \to G$ defined
by $\gamma(\Delta) = \{\beta^*_m(\pi_{U^m}(\Delta))\}$ (where $\Delta$ is a
generator of ${\check H}_n(M)$) is an isomorphism of
${\check H}_n(M)$ onto $G$.  We shall use the isomorphism in what follows
and for convenience we shall let $\gamma(\Delta) = \{z^n_m(\Delta)\}$, i.e\.
$z^n_m(\Delta) = \beta^*_m(\pi_{U^m}(\Delta)) \in H_n(V^m_n)$.  
\sk1

Let $\Delta$ be the generator of ${\check H}(M)$
where $\Delta = \{z^n_m(\Delta)\}$, a sequence of $n$-cycles such that $\pi^*_m : z^n_{m+1}(\Delta) \to z^n_m(\Delta)$ where  
$z^n_m(\Delta)$ is the $m^{\text{th}}$
coordinate of $\Delta$, i.e., $\pi_{V^m_n}(\Delta) = z^n_m(\Delta)$ an
$n$-cycle in $N(V^m_n)$ where $N(V^m_n)$ is the $n$-skeleton of $N(V^m)$.
If $s\Delta$, $s\in Z_p$, is any $n$-cycle in ${\check H}(M)$, then the
coordinate $n$-cycles $z^n_m(s\Delta)$ and $z^n_m(\Delta)$ contain exactly
the same $n$-simplices in $N(V^m_n)$ (see Lemma 8 below).

It follows from Lemma 2 that there is no loss of generality in assuming that 
$\pi_{V^m_n} : {\check H}(M) \to H_n(V^m_n)$ has the property that 
$\pi_{V^m_n} ({\check H}(M)) \cong Z_p$ for each natural number $m$.
\vskip 30pt

\cl{A Cernavskii Operator $\sigma$}
\sk1

\ni 4.  {\it An Operator $\sigma$ is defined on $n$-Chains similar to the
Cernavskii operator in [27].}
 
Choose $g\in A_p - H_1$.   
Next, define an operator $\sigma$ [cf\. 27, 10] on $n$-chains.
Recall that if $\delta^n$ is an $n$-simplex in $N(V^m)$, then either
$N[\delta^n]\subset$ interior of $F_\phi$ or $N[\delta^n] \cap F_\phi = 
\eset$.  If $\delta^n$ is an $n$-simplex in $N(V^{m+1})$ and int $F\supset 
N[\delta^n]$, then $\sigma_m(\delta^n) = \ds{\sum^{p-1}_{s=0}}g^s(\delta^n) 
= p\delta^n = 0$ mod $p$ since $g^s(\delta^n) = \delta^n$.
If $\delta^n$ is an $n$-simplex in $N(V^{m+1})$ and 
$N[\delta^n] \cap F = \eset$, then $\delta^n$ is in a unique non degenerate
distinguished family $S^{m+1}_i$ of $n$-simplices which has cardinality 
$p^{c_i}$.  To see this, observe that if $\delta^n$ is a singleton distinguished
family in $N(V^m)$, then the ``vertices'' of $\delta^n$ are members of $W^{m+1}_F$
and, consequently, $N[\delta^n] \cap F \neq \eset$ by construction of $W^{m+1}_F$.
Hence, $\delta^n$ is in a unique non degenerate distinguished family $S^{m+1}_i$ of
$n$-simplices.  If $\pi^*_m(\delta^n)$ is an $n$-simplex, then each $n$-simplex in
$S^{m+1}_i$ projects by $\pi^*_m$ to an $n$-simplex in $N(V^m)$.  If
$\pi^*_m(\delta^n)$ is a $k$-simplex, $k < n$, then each $n$-simplex in $S^{m+1}_i$
projects under $\pi^*_m$ to a $k$-simplex.  Let $\sigma
(\delta^n) = \ds{\sum^{p-1}_{s=0}}g^s(\delta^n)$ where $g^0$ is the identity
homeomorphism.  Observe that if $\delta^n$ is an $n$-simplex in $N(V^{m+1})$
such that $N[\delta^n] \cap F = \eset$ and $\pi^*_m(\delta^n)$ is an 
$n$-simplex in $N(V^m)$, then $N[\pi^*_m(\delta^n)] \cap F = \eset$ by the
construction of $W^{m+1}_F$ and either (1) $\pi^*_m(\delta^n)$ is a singleton family
in $N(V^m)$ in which case it is a singleton family consisting of a $k$-simplex with $k
< n$ by the construction of $W^m_F$ contrary to the assumption that $\pi^*_m(\delta^n)$
is an $n$-simplex or (2) $\pi^*_m(\delta^n)$ is in a non degenerate
distinguished family $S^m_j$ of $n$-simplices in $N(V^m)$. In Case (2), 
$\pi^*_m$ maps $\big\{g^s(\delta^n)\big\}^{p-1}_{s=0}$ one-to-one onto 
$\big\{\pi^*_m(g^s(\delta^n))\big\}^{p-1}_{s=0} = \big\{g^s(\pi^*_m
(\delta^n))\big\}^{p-1}_{s=0}$ and $\pi^*_m\sigma_m\delta^n = \sigma_m\pi^*_m
\delta^n$.  It should be clear that if $\delta^n$ is an $n$-simplex in
$N(V^{m+1})$, then $\pi^*_m\sigma_m\delta^n = \pi^*_m\ds{\sum^{p-1}_{s=0}}
g^s(\delta^n) = \ds{\sum^{p-1}_{s=0}}\pi^*_mg^s(\delta^n) = 
\ds{\sum^{p-1}_{s=0}}g^s(\pi^*_m(\delta^n)) = \sigma_m\pi^*_m\delta^n$.  
Clearly, $\sigma_m$ can be extended to any $n$-chain in $N(V^m_n)$.  If 
$\delta^n$ is an $n$-simplex in $N(V^m_n)$ and $\pi^*_m\delta^n$ is a 
$k$-simplex, $k < n$, then $\pi^*_m\sigma_m(\delta^n)$ is a trivial 
$n$-chain.  

Recall that the members of a
distinguished family of $k$-simplices in $N(V^m)$
have the same orientation (being the lifting of a
$k$-simplex in $N(R^m)$ where $R^m$ is a certain cover of $Y$).  

Recall that if $\delta^n$ is an $n$-simplex in $N(V^m)$, then either
interior $F\supset N[\delta^n]$ or $N[\delta^n] \cap F = \eset$.  This
follows from the construction of $V^m$.

\ni Lemma 6. {\rm{[cf\. 27]}}  
The special operator $\sigma_m$ 
maps $n$-cycles to $n$-cycles and $\sigma_m$ commutes with the special
projections on $n$-chains.  If $z^n_{m+1}(\Delta) = z$ is a coordinate 
$n$-cycle in $N(V^{m+1}_n)$, then either (a) $\sigma_mz = 0$ or (b) $\sigma(z)
= z$.  

\ni Proof.  The homeomorphism $g$ induces a simplicial homeomorphism 
$\bar g$ of $N(V^{m+1})$ onto itself since $g$ maps ``vertices'' (elements
of $V^{m+1})$ one-to-one onto ``vertices''.  

Let $\Delta$ be a non zero $n$-cycle in ${\check H}_n(M) =
\varprojlim H_n(V^{m+1}_n) \cong Z_p$.  Let $z = z^n_{m+1}(\Delta)$ be the
coordinate $n$-cycle of $\Delta$ in $H_n(V^{m+1}_n)$.  Let $\sigma =
\sigma_m$.  Consider $\sigma(z)$ where $z = \ds{\sum^k_{i=1}}c_i\delta^n_i$.
Thus, $\sigma z = \ds{\sum^k_{i=1}}c_i\sigma \delta^n_i$  where 
$\sigma\delta^n_i = \ds{\sum^{p-1}_{s=0}}g^s(\delta^n_i)$ and
$\big\{g^s(\delta^n_i)\big\}^{p-1}_{s=0}$ is the distinguished subfamily of
$n$-simplices in $N(V^{m+1})$ associated with $\delta^n_i$.
It follows
by the construction that there is $z^n_{m+2}(\Delta)$, a coordinate
$n$-cycle in $N(V^{m+2}_n)$ such that $g^s(z^n_{m+2}(\Delta))$ maps by
$\pi^*_{m+1}$ onto $g^s(z)$.  Hence, $g^s(z)$ for $0\leq s\leq p-1$ contains
the same $n$-simplices as $z$ (see Lemma 8).  
Hence, $\sigma z = \ds{\sum^{p-1}_{t=0}}z_t$ where $z_0 = z$ and $z_t = 
\ds{\sum^k_{i=1}}c_ig^t(\delta^n_i)$ where $g^t(\delta^n_i)$ means $\delta^n_i$
whenever int $F\supset N[\delta^n_i]$.  Since $g$ is a homeomorphism on $M$,
$g$ induces an automorphism on $\pi_{V^{m+1}}(G)$ in $H_n(V^m_n)$ where 
$G = \varprojlim H_n(V^m_n) \cong Z_p$, $\pi_{V^{m+1}_n}(G) \cong Z_p$, and
$z_t$ is an $n$-cycle.  It follows that $\sigma z = \ds{\sum^{p-1}_{s=0}}
g^s(z)$ which is an $n$-cycle.  If $g$ induces the identity automorphism, then
$\sigma z = pz = 0$ mod $p$ (the trivial $n$-cycle).  If the induced 
automorphism is not the identity, then it will be shown that 
$\ds{\sum^{p-1}_{t=1}} z_t = 0$ mod $p$ and that $\sigma z = z$, that is,
$\sigma$ is the identity automorphism.  Now, $\pi_{V^{m+1}_n}(G) \cong Z_p
= \{0,1,2,\cdots,p-1\}$.  Since $g$ induces an automorphism on 
$\pi_{V^{m+1}_n}(G)$, it induces an automorphism $g_*$ on $Z_p$. Let $g_*(1)
= x$. Hence, $g^s_*(1) = x^s$.  It is well known in number theory that $p$
divides $x^{p-1}-1$.  Also, $x^{p-1}-1 = (x-1)(1+x+x^2+\cdots + x^{p-2})$.
Consequently, $p$ divides $x(1+x+x^2+\cdots+x^{p-2}) = x + x^2 + \cdots + 
x^{p-1} = \ds{\sum^{p-1}_{s=1}}x^s = 0$ mod $p$.  It follows that $\sigma(z)
= \ds{\sum^{p-1}_{s=0}}g^s(z) = z + \ds{\sum^{p-1}_{s=1}}g^s(z) = z$ mod $p$.  
Thus, if $g$ does not induce the identity automorphism, then $\sigma$ is the
identity automorphism.   

If under the projection $\pi^*_m : H_n(V^{m+1}_n) \to H_n(V^m_n)$, the image 
of an $n$-simplex $\delta^n$ such that $N[\delta^n] \subset M - F$ 
is a $k$-simplex with $k < n$, then the same is true for all members of the
distinguished family to which $\delta^n$ belongs.  Thus,
$\pi^*_m\sigma_m\delta^n = \sigma_m\pi^*_m\delta^n = 0$ where $\sigma_m
\pi^*_m\delta^n$ is defined above and is a trivial $n$-chain.  
If $\pi^*_m(\delta^n)$ is an $n$-simplex, then by construction the
distinguished subfamily of $n$-simplices with 
which $\delta^n$ is associated is in one-to-one
correspondence with the distinguished subfamily of $n$-simplices
with which $\pi^*_m(\delta^n)$
is associated.  Thus, $\pi^*_m\sigma_m\delta^n = \sigma_m\pi^*_m\delta^n$.  
If $N[\delta^n] \subset F_\phi$, then
$\pi^*_m\sigma_m\delta^n = 0 = \sigma_m\pi^*_m\delta^n$.  Consequently, $\sigma$
carries over to the $n$-cycles of $M$ and to ${\check H}(M)$.  Lemma 6 is 
proved.

\ni Lemma 7.  If $\Delta_1$ and $\Delta_2$ are non zero elements of ${\check
H}_n(M)$, then for each $m$, exactly the same simplices appear in the chains
$z^n_m(\Delta_1)$ and $z^n_m(\Delta_2)$ which are in the $n$-dimensional
complex, $N(V^m_n)$, the $n$-skeleton of $N(V^m)$.

\ni Proof.  Suppose that this is not true.  Without loss of generality,
assume that the $n$-simplex $\delta^n$ appears in $z^n_m(\Delta_1)$ and not
in $z^n_m(\Delta_2)$ where $z^n_m(\Delta_1)$ and $z^n_m(\Delta_2)$ are the
$m^{\text{th}}$ coordinates of $\Delta_1$ and $\Delta_2$, respectively.
Since $N(V^m_n)$ is $n$-dimensional, these coordinates are $n$-cycles in
$N(V^m_n)$.  Since ${\check H}_n(M) = Z_p$, assume that $\Delta_1$
generates ${\check H}_n(M)$ and that $\Delta_2 = s\Delta_1$ for some natural
number $s$, $1\leq s < p$.  
Consequently, $z^n_m(\Delta_2) = s(z^n_m(\Delta_1))$.
It follows that $\delta^n$ appears in $s(z^n_m(\Delta_1))$ and hence in
$z^n_m(\Delta_2)$ -- a contradiction.

\ni Lemma 8.  Suppose that $\Delta \in {\check H}_n(M)$ with $\Delta \neq 0$
and $z = z^n_{m+1}(\Delta) = \pi_{V^{m+1}_n}(\Delta)$, the coordinate 
$n$-cycle of $\Delta$ in $H_n(V^{m+1}_n)$.  Let $z =\ds{\sum^q_{i=1}}c_i
\delta^n_i$.  Then the collection $C_j = \{\delta^n_{j_1},\delta^n_{j_2},
\cdots,\delta^n_{j_{t_j}}\}$ of all $n$-simplices in $\{\delta^n_1,\delta^n_2,
\cdots,\delta^n_q\}$ which are in a fixed distinguished family $S^{m+1}_j$
of $n$-simplices in $N(V^{m+1}_n)$ have the properties (1) if $x = 
\ds{\sum^t_{i=1}}c_{j_i}\delta^n_{j_i}$, then either (a) $\sigma_mx = 0$ when 
$g$ induces the identity automorphism or (b) $\sigma_mx = x$ when $g$ does not
induce the identity automorphism and (2) 
$\ds{\sum^t_{i=1}}c_{j_i} = 0$ mod $p$.

\ni Proof.  By Lemma 6, either (a) $\sigma_mz = 0$ or (b) $\sigma_m(z) = z$.

Case (a):  Since $\sigma_m(z) = 0$, it follows that $\sigma_mx = 0$
since for each $i$, $1\leq i \leq t$, $\sigma_m \delta^n_{j_i} = 
\ds{\sum^{p-1}_{s=0}}g^s(\delta^n_{j_i})$ where for each $s$, $0\leq s \leq
p-1$, $g^s(\delta^n_{j_i})$ is an $n$-simplex in $C_q\subset S^{m+1}_j$.  
Choose notation such that (1) $S^{m+1}_i = \big\{\delta_j\big\}^{p^k}_{j=1}$
($p^k = p^{c_j}$ in earlier notation), (2) $g\delta_j =\delta_{j+1}$ for 
$1\leq j < p^k$ and $g\delta_{p^k} = \delta_1$ ($g$ permutes the $\delta_j$ 
in a cyclic order), and (3) $x = \ds{\sum^{p^k}_{i=1}}c_i\delta_i$ where
$c_i = 0$ iff $\delta_i \not\in C_q$ and $c_i = c_{j_i}$ iff $\delta_i =
\delta^n_{j_i} \in C_q$.  Let $\sigma = \sigma_m$.  Since $\sigma(z) = 0$, it
follows that $\sigma(x) = 0$.  Note that $0 = \sigma(x) = c_1\ds{\sum^p_{i=1}}
\delta_i + c_2 \ds{\sum^p_{i=1}}\delta_{i+1} + \cdots + c_j \ds{\sum^p_{i=1}}
\delta_{i+j-1} + \cdots + c_{p^k-p} \ds{\sum^p_{i=1}}\delta_{i+p^k-p-1} + 
\cdots + c_{p^k}(\delta_{p^k} +\delta_1 + \delta_2 + \cdots + \delta_{p-1})$.
Rearrange as 
$$\eqalign
{\sigma(x) &= (c_1 + c_{p^k} + c_{p^k-1} + \cdots + c_{p^k-p+2})
\delta_1 + \cr
&(c_2 + c_1 + c_{p^k} + \cdots + c_{p^k-p+1})\delta_2 + \cr
&\vdots\cr
&(c_p + c_{p-1} + c_{p-2} + \cdots + c_1)\delta_p +  \cr
&(c_{p+1} + c_p + c_{p-1} + \cdots + c_2)\delta_{p+1} + \cr
&\vdots\cr
&(c_{p^k} + c_{p^k-1} + \cdots + c_{p^k-p+1})\delta_{p^k}.\cr}$$
Since $\sigma(x) = 0$, it follows that the coefficient of $\delta_i$ is 0
mod $p$ for $1\leq i \leq p^k$.  A careful consideration of pairs of successive
coefficients of $\delta_i$ and $\delta_{i+1}$ will give the following result.
If $1\leq i \leq p^k$, $1\leq j \leq p^k$, and $i \equiv j$ mod $p$, then
$c_i \equiv c_j$ mod $p$. Thus $\ds{\sum^{p^k}_{i=1}}c_i = \ds{\sum^p_{i=1}}
c_i + \ds{\sum^{2p}_{i=p}}c_i + \ds{\sum^{3p}_{i=2p+1}}c_i + \cdots + 
\ds{\sum^{(p^{k-1})p}_{i=(p^{k-1}-1)p}}c_i$ with $p^{k-1}$ summations each of 
length $p$. Now, $c_i\equiv c_j$ mod $p$ if $i\equiv j$ mod $p$ gives that each
of the $p^{k-1}$ summations is congruent to 0 mod $p$.  Thus, 
$\ds{\sum^{pk}_{i=1}} c_i = p^{k-1}\left(\ds{\sum^p_{i=1}}c_i\right) \equiv 0$
mod $p$ if $k > 1$.
If $k = 1$, then $x = \ds{\sum^p_{i=1}}c_i\delta_i$, 
$$\eqalign
{\sigma(x) &= c_1(\delta_1 + \delta_2 + \cdots + \delta_p)\cr
&+ c_2(\delta_2+\delta_3 + \cdots + \delta_{p+1})\cr
&\vdots\cr
&+ c_p(\delta_p + \delta_1+\delta_2 +\cdots + \delta_{p-1}),\ 
\hbox{ rearrange as}
\cr
&(c_1 + c_2 + \cdots + c_p)\delta_1 +\cr
&(c_1 + c_2 + \cdots + c_p)\delta_2 + \cr
&\vdots\cr
&(c_1+c_2 + \cdots + c_p)\delta_p,\cr}$$
and $\ds{\sum^p_{i=1}}c_i\equiv 0$ mod $p$.

It is instructive to consider a simple example.  Let $p=3$ and $x = 
\ds{\sum^9_{i=1}}c_i\delta_i$.  Thus, $\sigma x = c_1(\delta_1+\delta_2+
\delta_3) + c_2(\delta_2 + \delta_3 + \delta_4) + c_3(\delta_3+\delta_4
+\delta_5) + 
c_4(\delta_4 + \delta_5 + \delta_6) + c_5 (\delta_5 + \delta_6 + \delta_7)
+ c_6(\delta_6 + \delta_7 + \delta_8) + c_7(\delta_7 + \delta_9 +\delta_9) + 
c_8(c_8 + c_9 + c_1) + c_9(c_9 + \delta_1+\delta_2) =$ (by rearrangement) $=
(c_1 + c_8 + c_9)\delta_1 + (c_1 + c_2 + c_9)\delta_2 + (c_1 + c_2 + c_3)
\delta_3 + (c_2+c_3 + c_4)\delta_4 + (c_3 + c_4 + c_5)\delta_5 + (c_4 + c_5 +
c_6)\delta_6 + (c_5 + c_6 + c_7)\delta_7 + (c_6+c_7 + c_8)\delta_8 + (c_7 + c_8 + c_9\delta_9$.  For each $i$, $1\leq i \leq 9$, the coefficient of $\delta_i
= 0$ mod 3.  Observe that from the coefficients of $\delta_1$ and $\delta_2$,
it follows that $c_2 \equiv c_8$ mod 9.  The coefficients of $\delta_2$ and
$\delta_3$ yield that $c_3\equiv c_9$ mod 3.  Continuing, $c_1\equiv c_4$,
$c_2 \equiv c_5$, $c_3 \equiv c_6$, $c_4 \equiv c_7$, $c_5 \equiv c_8$,
$c_6 \equiv c_9$, and $c_7 \equiv c_1$ all mod 3.  Thus, $(c_1 + c_2 + c_3)
+ (c_4 + c_5 + c_6) + (c_7 + c_8 + c_9) = 0$ mod 3 since $(c_4 + c_5 + c_6)
\equiv (c_1+c_2+c_3)$ mod 3, $(c_7 + c_8 + c_9) \equiv (c_4 + c_5 + c_6) \equiv
(c_1+c_2+c_3)$ mod 3 and $\ds{\sum^9_{i=1}}c_i \equiv 3(c_1+c_2+c_3)
\equiv 0$ mod 3.

Case (b): $\sigma_m(x) = x$.  Choose notation as in Case (a).  Write $\sigma
(x)$ as in Case (a), but in this case, $x = \sigma(x)$ rather than $0 = 
\sigma(x)$.  Consider first the example $p=3$ and $x = \ds{\sum^9_{i=1}}
c_i\delta_i$ where $k = 2$.  Now, $\sigma x = (c_1 + c_8 + c_9)\delta_1 + 
(c_1+c_2+c_9)\delta_2 + (c_1+c_2+c_3)\delta_3 + (c_2 + c_3 + c_4)\delta_4 + 
(c_3 + c_4 + c_5)\delta_5 + (c_4 + c_5 + c_6)\delta_6 + (c_5 + c_6 + c_7)
\delta_7 + (c_6 + c_7 + c_8)\delta_8 + (c_7 + c_8 + c_9) \delta_9 = x = 
c_1 \delta_1 + c_2\delta_2 + c_3\delta_3 + c_4 + \delta_4 + c_5 \delta_5
+ c_6 \delta_6 + c_7\delta_7 + c_8\delta_8 + c_9\delta_9$.  This is an 
identity.  Thus, the coefficient of $\delta_i$ on one side is equal mod $p$
to the coefficient of $\delta_i$ on the other side.  Hence, $c_1 + c_8 + c_9
\equiv c_1$ mod $p$, $c_1 + c_2 + c_9 \equiv c_2$ mod $p$, $c_1 + c_2 + c_3
\equiv c_3$ mod $p$, and so forth.  Thus, $\ds{\sum^9_{i=1}}c_i\equiv (c_1
+ c_8 + c_9) + (c_1 + c_2 + c_9) + (c_1+c_2+c_3) + (c_2 + c_3 + c_4) + 
(c_3 + c_4 + c_5) + (c_4 + c_5 + c_6) + (c_5 + c_6 + c_7) + (c_6 + c_7 + c_8)
+ (c_7 + c_8 + c_9) = 3(c_1+c_2+c_3+c_4+c_5+c_6+c_7+c_9) \equiv 0$ mod $p$.  

Consider the general case as in Case (a) but with $x = \sigma(x)$ rather than
$0 = \sigma(x)$.  Hence, 
$$\eqalign
{x = \sigma(x) &= (c_1 + c_{p^k} + c_{p^k-1} + \cdots + c_{p^k-p+2})\delta_1
+\cr
&(c_2+c_1 + c_{p^k} + \cdots + c_{p^k-p+1}\delta_2 +\cr
&\vdots\cr
&(c_p + c_{p-1} + c_{p-2} + \cdots + c_1)\delta_p +\cr
&(c_{p+1} + c_p + c_{p-1} + \cdots + c_2\delta_{p+1} +\cr
&\vdots\cr
&(c_{p^k} + c_{p^k-1} + \cdots + c_{p^k-p+1})\delta_{p^k} = \sum^{p^k}_{i=1}
c_i\delta_i.\cr}$$
It follows from this identity that the coefficient of $\delta_i$ on one side
is equal mod $p$ to the coefficient of $\delta_i$ on the other side.  
Consequently, $\ds{\sum^{p^k}_{i=1}}c_i = p \ds{\sum^{p^k}_{i=1}}c_i$ mod
$p = 0$ mod $p$ as claimed where $k > 1$. For $k = 1$, $x = \ds{\sum^p_{i=1}}
c_i\delta_i = \left(\ds{\sum^p_{i=1}}c_i\right) \delta_1 + 
\left(\ds{\sum^p_{i=1}}c_i\right) \delta_2 + \cdots + \left(\ds{\sum^p_{i=1}}
c_i\right) \delta_p$, $c_t \equiv \ds{\sum^p_{i=1}}c_i$ mod $p$ for each $t$,
$1\leq t\leq p$, and $\ds{\sum^p_{i=1}}c_i \equiv p\left(\ds{\sum^p_{i=1}}
c_i\right)$ mod $p = 0$ mod $p$.  Lemma 8 is proved.

\ni Lemma 9.  {\rm{[cf\. 27]}}  If $W$ is any domain in $M$ (connected open 
set), then ${\check H}_n(M-W) = 0$ and if $\Delta \in {\check H}_n(M)$ with
$\Delta \neq 0$, then there is a natural number $m_0$ such that for each
$m > m_0$, the carrier of $z_m(\Delta)$ meets $W$, that is, if $z_m =
\ds{\sum^r_{j=1}}g_j\delta^n_j$ where $1\leq g_j < p$, then for some $j$,
$N[\delta^n_j] \cap W \neq \eset$ where $N[\delta^n_j]$ is the nucleus
of $\delta^n_j$.  

\ni Proof.  Here, the fact that ${\check H}_n(M-W) = 0$ is used.  Suppose
that the lemma is false.  There is a subsequence $\{V^{m_i}\}$ cofinal in the
collection of all open coverings of $M$ such that $z_{m_i}(\Delta)$ is an
$n$-cycle in $N(V^{m_i}\mid (M-W))$.  Thus, $z_{m_i}(\Delta)$ is not a 
bounding chain in $N(V^{m_i})$ and, hence not a bounding chain in $N(V^{m_i}
\mid (M-W))$.  Thus, $z_{m_i}(\Delta)$ determines a non zero element in
$H_n(V^{m_i}\mid (M-W))$ and $\{z_{m_i}(\Delta)\}$ determines a non zero 
element of ${\check H}_n(M-W)$.  This yields a contradiction.  The lemma is
proved.

\ni Lemma 10. {\rm{[cf\. 27]}}  
If $K =$ interior $F_\phi$, then $K = \eset$ and $F_\phi$ is
nowhere dense.

\ni Proof.  Suppose that $K\neq \eset$.  Let $\Delta$ be a non zero
$n$-cycle in ${\check H}_n(M) \cong Z_p$ and $z^n_m(\Delta) =
\beta^*_m\pi_{U^m}(\Delta) \in H_n(V^m_n)$.  Let $z^n_m(\Delta) = C^m_1 +
C^m_2$ where $C^m_1$ is an $n$-chain such that each $n$-simplex in $C^m_1$
has a nucleus in $M - \bar K$ and $C^m_2$ is an $n$-chain such that each
$n$-simplex in $C^m_2$ has a nucleus in $K$.  It will be shown that both
$C^m_1$ and $C^m_2$ are $n$-cycles.  Let $z = z^n_m(\Delta)$.

\ni Proof.  It was shown in the proof of Lemma 6 that either (a)
$\sigma z = 0$ or (b) $\sigma z = z$.  Let
$z = \ds{\sum^k_{i=1}}c_i\delta^n_i$.  
If int $F \supset N[\delta^n_i]$, then $g^s(\delta^n_i) = \delta^n_i$ for
each $s$.  Suppose that $\delta^n_j$ shares an $(n-1)$ face $\delta^{n-1}$
with $\delta^n_i$ where int $F\not\supset N[\delta^n_i]$.  Then $g^s(\delta^n_j)
\neq g^t(\delta^n_j)$ for $s \neq t$ but $g_s(\delta^n_j)$ and
$g^t(\delta^n_j)$ share the $(n-1)$-face $\delta^{n-1}$.  Let $C_j =
\{\delta^n_{j_1},\delta^n_{j_2},\cdots,\delta^n_{j_t}\}$ denote the 
collection of all $n$-simplices which (1) appear in $z$ (with non zero
coefficient) where for each $c$, $1\leq c \leq t$, $\delta^n_{j_c}$ shares
the $(n-1)$-face $\delta^{n-1}$ with $\delta^n_i$ and (2) have the property
that the non degenerate distinguished family $S^m_j$ of $n$-simplices in 
$N(V^m)$ contains $C_j$.  By Lemma 8, $\ds{\sum^t_{c=1}}c_{j_c} = 0$ mod $p$.
Thus, the coefficient of $\delta^{n-1}$ in $\partial C^m_1$ is 0 mod $p$. 

To show that $\partial C^m_1$ is 0, it suffices to show that each
such $(n-1)$ simplex $\delta^{n-1}$ in $\partial C^m_1$ which is a face of
some $\delta^n_i$ where int $F\supset N[\delta^n_i]$ is 0.  As shown above,
this is the case and $C^m_1$ is an $n$-cycle.  
Thus, $\partial C^m_1 = 0$ and $C^m_1$ is an $n$-cycle.  It follows that
$C^m_2$ is an $n$-cycle.  
From the definition of the special projections, $\pi^*_m(C^{m+1}_i)
= C^m_i$ for $i = 1,2$.  Thus, we can write $\Delta = \Delta_1 + \Delta_2$
where $z^n_m (\Delta_i) = C^m_i$ for $i = 1,2$.  Now, the nucleus of each
simplex in $z^n_m(\Delta_2)$ misses the nonempty open set $M - \bar K$ for
each $m$, so by Lemma 7, $\Delta_2 = 0$.  Similarly, the nucleus of each
simplex in $z^n_m(\Delta_1)$ misses the non-empty open set $K$ for each $m$.
Hence, by Lemma 7, $\Delta_1 = 0$. Consequently, we have $\Delta = 0$, a
contradiction. The lemma is proved.   
\vskip 30pt

\ni 5.  {\it A proof that a $p$-adic group $A_p$ can not act effectively 
on a compact connected $n$-manifold where $\phi : M \to M/A_p$ is the
orbit mapping.}

{\bf Remarks.}  If the compact connected $n$-manifold $M$ has a non empty
boundary, then two copies of $M$ can be sewed together by identifying the
boundaries in such a way that the result is a compact connected 
$n$-manifold $M'$ without boundary.  If $A_p$ acts effectively on
$M$, then $A_p$ acts effectively on $M'$.  
If $M$ is not an orientable $n$-manifold, then we can take the
double cover of $M$ on which $A_p$ acts effectively if it acts effectively on $M$.
There is no loss of generality in assuming that $M$ is a compact connected
orientable $n$-manifold without boundary.

\ni {\bf Definition.}  An $n$-manifold $(M,d)$ is said to have Newman's
Property w.r.t\. the class $L(M,p)$ (as stated above) iff 
there is $\epsilon > 0$ such that for any $\phi \in L(M;p)$, there is 
some $x\in
M$ such that diam $\phi^{-1}\phi(x) \geq \epsilon$.

Generalizations can be made to metric spaces $(X,d)$ which are locally
compact, connected, and $lc^n$ [4] which have domains $D$ such that $\bar
D$ is compact, $lc^n$, and $H_n(X,X-D),Z_p) \cong Z_p$.  

\ni Theorem.  If $L(M,p)$ is the class of all orbit mappings $\phi : M \to
M/A_p$ where $A_p$ acts effectively on a compact, connected, and orientable
$n$-manifold $M$, then $M$ has Newman's Property w.r.t\. $L(M,p)$.  

\ni Proof.  There is no loss of generality in assuming that $M$ is orientable
and has empty boundary.  

By hypothesis, $\check H_n(M) \cong Z_p$.  
Consider a finite open covering $U = W_1$ where $W_1$ and
$W_2$ satisfy Lemma 2 and such that if $z(\Delta)$ is the $V$-coordinate of
a non-zero $n$-cycle $\Delta \in \check H_n(M)$ where $V$ refines $W_2$, then
$\pi_{VU} z(\Delta) \neq 0$.  Let $\epsilon$ be the Lebesque number of
$W_2$.  Suppose that there is $\phi \in L(M,p)$ such that diam
$\phi^{-1}\phi(x) < \epsilon$ for each $x\in M$.  Construct the special
coverings $\{V^m\}$ and the special refinements $\{U^m\}$ as in Lemma 5 such
that the star of each distinguished family of $V^1$ lies in some element of
$W_2$.  Furthermore, the special projections $\pi_m$ can be constructed such
that if $\big\{\delta^n_{sj}\big\}^{t^m_s}_{j=1}$ is a distinguished family
of $n$-simplices in $N(V^m_n)$, then $\pi_{V^mU}$ takes $\delta^n_{sj}$,
$1\leq j \leq t_s$, to the same simplex $\delta_s$ in $N(U)$.  
Now, let $z_m = z_m^n(\Delta)$.  

Let $z_m = \ds{\sum^k_{i=1}}c_i\delta^n_i$.  By Lemma 10, $F\cap
N[\delta^n_i] = \eset$ for each $i$ since int $F = \eset$.  Hence, for each
$j$, $1\leq j \leq k$, $\delta^n_j$ is in a non degenerate
distinguished family $S^m_j$
of $n$-simplices in $N(V^m)$.  Let $C_j = \{\delta^n_{j_1},\delta^n_{j_2},
\cdots,\delta^n_{j_t}\}$ denote the collection of all $n$-simplices such that
(1) $\delta^n_{j_i}$ appears in $z_m$ for $1\leq i \leq t$ and (2) $S^m_j
\supset C_j$.  By Lemma 8, $\ds{\sum^t_{i=1}}c_{j_i} = 0$ mod $p$.  Since
the $n$-simplices in $C_j$ are sent by $\pi_{V^mU}$ to a single simplex 
$\delta_j$ in $N(U)$, it follows that the coefficient of $\delta_j$ is 0 mod
$p$ and $z$ is sent by $\pi_{V^mU}$ to the zero $n$-cycle in $N(U)$.  
Thus, the projection of $z_m$ by 
$\pi_{V^mU} : H_n(V^m_n) \to H_n(U)$
takes the nontrivial $n$-cycle $z_m(\Delta)$ to the 0 $n$-cycle mod $p$.  This
violates the conclusion of Lemma 2.  Thus, $M$ has Newman's Property w.r.t\. the
class $L(M,p)$.  Hence, $\epsilon$ is a Newman's number and the Theorem is
proved.  

It is well known that 
if $A_p$ acts effectively on a compact connected $n$-manifold
$M$, then given any $\epsilon > 0$, there is an effective action of $A_p$
on $M$ such that diam $\phi^{-1}\phi(x) < \epsilon$ for each $x\in M$.  That
is, $M$ fails to have Newman's property w.r.t. $L(M,p)$.  It follows that
$A_p$ can not act effectively on a compact connected $n$-manifold $M$.  
\vskip 30pt

\ni 6.  {\it How to obtain a proof that a $p$-adic group can not act effectively on
a connected $n$-manifold}.

As indicated above, there is no loss of generality in assuming that $M$ is a
connected orientable $n$-manifold without boundary. If $A_p$ acts effectively on
$M$ (which is locally compact), then the orbit map $\phi : M\to M/A_p$ is open and
closed with $\phi^{-1}\phi(x)$ compact for each $x\in M$.  Hence, $\phi$ is a
proper map (if $M/A_p\supset A$ and $A$ is compact, then $\phi^{-1}(A)$ is
compact).  

Construct sequences $\{V^m\}$ and $\{U^m\}$ of locally finite open coverings of $M$
by constructing locally finite open coverings $R^m$ of $M/A_p$ in the same manner
as in Lemma 3, 4, and 5 where each $r\in R^m$ has a compact closure.  Since $\phi$
is proper, the distinguished families of open sets in $V^m$ generated by members of
$R^m$ have the same properties as in Lemmas 1, 3, 4, and 5.  The proof follows as
in the compact case.  

Consequently, the {\it Hilbert Smith Conjecture} is true.

\newpage

References

\ni 1. Bochner, S. and Montgomery, D. Locally compact groups
of differentiable transformations,  Annals of Math 47 (1946)
639-653.
\sk1

\ni 2. Bredon, Glen E., Orientation In Generalized Manifolds And
Applications To The Theory of Transformation Groups, Mich. Math. Jour.
7 (1960) 35-64.
\sk1     

\ni 3. Browder, F., editor Mathematical development arising
from Hilbert Problems, in Proceedings of Symposia in Pure Mathematics
Northern Ill. Univ., 1974, XXVIII, Parts I and II.
\sk1

no 4. Floyd, E.E., Closed coverings in {\v C}ech homology
theory, 
TAMS 84 (1957) 319-337.
\sk1

\ni 5. Gleason, A., Groups without small subgroups,  Annals
of Math. 56 (1952) 193-212.
\sk1

\ni 6. Greenberg, M. Lectures On Algebriac Topology, 
New York, 1967.
\sk1

\ni 7.
Hall, D.W. and Spencer II, G.L.,
Elementary Topology,
John Wiley \& Sons, 1955.
\sk1

\ni 8. Hilbert, David, Mathematical Problems, BAMS 8
(1901-02) 437-479.
\sk1

\ni 9. Lickorish, W.B.R., The piecewise linear unknotting of
cones, Topology 4 (1965) 67-91.
\sk1

\ni 10. 
McAuley, L.F. and Robinson, E.E., 
On Newman's Theorem
For Finite-to-one Open Mappings on Manifolds, PAMS 87 (1983)
561-566.
\sk1

\ni 11.
McAuley, L.F. and Robinson, E.E., 
Discrete open and closed mappings on generalized
continua and Newman's Property,  Can. Jour. Math. XXXVI (1984)
1081-1112.
\sk1

\ni 12. Michael, E.A., Continuous Selections II, Annals of
Math. 64 (1956) 562-580.
\sk1

\ni 13. Montgomery, D. and Zippin, L., 
Small groups of
finite-dimensional groups, Annals of Math. 56 (1952) 
213-241.
\sk1

\ni 14. Montgomery, D. and Zippin, L., 
Topological Transformation Groups, Wiley
(Interscience), N.Y., 1955.
\sk1

\ni 15. Newman, M.H.A., A theorem on periodic transformations of
spaces, Q. Jour. Math. 2 (1931) 1-9.
\sk1

\ni 16. Robinson, Eric E., A characterization of certain
branched coverings as group actions, Fund. Math. CIII (1979)
43-45.
\sk1

\ni 17. Smith, P.A., Transformations of finite period. III
Newman's Theorem, Annals of Math. 42 (1941) 446-457.
\sk1

\ni 18. Walsh, John J., Light open and open mappings on
manifolds II, TAMS 217 (1976) 271-284.
\sk1

\ni 19. Whyburn, G.T., Topological Analysis,
Princeton University Press, Princeton, N.J., 1958.
\sk1

\ni 20. Wilson, David C., Monotone open and light open
dimension raising mappings, Ph.D. Thesis, Rutgers University,
(1969). 
\sk1

\ni 21. Wilson, David C., Open mappings on manifolds and a counterexample to
the Whyburn Conjecture Duke Math. Jour. 40 (1973) 705-716.
\sk1

\ni 22. Yang, C.T., Hilbert's Fifth Problem and related problems
on transformation groups, in Proceedings of Symposia in Pure Mathematics, 
Mathematical development arising from Hilbert Problems
Northern Illinois University, 1974, edited by F. Browder,
XXVIII, Part I (1976) 142-164. 
\sk1

\ni 23. Yang, C.T., $p$-adic transformation groups, Mich. Math. J.
 7 (1960) 201-218.
\vskip 15pt

\cl {Some Additional References To Work Associated With Efforts Made}
\cl{Towards Solving Hilbert's Fifth Problem (General Version)}
\sk1

\ni 24. Anderson, R.D.,  Zero-dimensional compact groups of
homeomorphisms,  Pac. Jour. Math. 7 (1957) 797-810.
\sk1

\ni 25.  Bredon, G.E., Raymond, F., and Williams, R.F.,  $p$-adic
groups of Transformations, TAMS 99 (1961) 488-498.
\sk1

\ni 26. Brouwer, L.E.J., 
Die Theorie der endlichen kontinuierlichen Gruppen unabhagin von
Axiom von Lie, Math. Ann.  67 (1909)
246-267 and 69 (1910), 181-203.
\sk1

\ni 27. Cernavskii, A.V.,
Finite-to-one open mappings of manifolds, Trans. of Math. Sk.
65 (107) (1964).
\sk1

\ni 28. Coppola, Alan J., On $p$-adic transformation groups,
Ph.D. Thesis, SUNY-Binghamton, 1980.
\sk1

\ni 29. Kolmogoroff, A. {\" U}ber offene Abbildungen, Ann.
of Math. 2(38) (1937) 36-38.
\sk1

\ni 30. Ku, H-T, Ku, M-C, and Mann, L.N.,  Newman's Theorem And
The Hilbert-Smith Conjecture, Cont. Math. 36 (1985) 489-497.
\sk1

\ni 31. Lee, C.N., Compact $0$-dimensional transformation group, 
(cited in [23] prior to publication, perhaps, unpublished).
\sk1

\ni 32. McAuley, L.F., Dyadic Coverings and Dyadic Actions, 
Hous. Jour. Math. 3 (1977) 239-246.
\sk1

\ni 33. McAuley, L.F., $p$-adic Polyhedra and $p$-adic Actions, 
Topology Proc. Auburn Univ. 1 (1976) 11-16.
\sk1

\ni  34. Montgomery, D., 
Topological groups of differentiable transformations, 
Annals of Math.  46 (1945) 382-387.
\sk1

\ni 35. Pontryagin, L.,  Topological Groups, (trans. by E. Lechner), 
Princeton Math. Series Vol. 2, Princeton University Press, Princeton, N.J.
(1939).
\sk1

\ni 36. Raymond, F., The orbit spaces of totally disconnected
groups of transformations on manifolds, PAMS 12 (1961) 1-7.
\sk1

\ni 37. Raymond, F.,  Cohomological and dimension theoretical properties
of orbit spaces of $p$-adic actions, in Proc. Conf. Transformation
Groups, New Orleans, (1967) 354-365.
\sk1

\ni 38. Raymond, F.,  Two problems in the theory of generalized
manifolds, Mich. Jour Math. 14 (1967) 353-356.  
\sk1

\ni 39. Raymond, F. and Williams, R.F.,  Examples of $p$-adic
transformation groups, Annals of Math. 78 (1963) 92-106.
\sk1

\ni 40. von Neumann, J., 
Die Einfuhrung analytischer Parameter in topologischen Gruppen,
Annals of Math.  34 (1933) 170-190.
\sk1

\ni 41. Williams, R.F.,  An useful functor and three famous
examples in Topology, TAMS 106 (1963) 319-329.
\sk1

\ni 42. Williams, R.F., The construction of certain $0$-dimensional
transformation groups, TAMS 129 (1967) 140-156.
\sk1

\ni 43. Yamabe, H., 
On a conjecture of Iwasawa and Gleason, Annals of Math.
58 (1953) 48-54.
\sk1

\ni 44. Yang, C.T., Transformation groups on a homological
manifold, TAMS 87 (1958) 261-283.
\sk1

\ni 45. Yang, C.T., $p$-adic transformation groups, Mich. Math.
Jour. 7 (1960) 201-218.
\sk1

\ni 46. Zippin, Leo, Transformation groups,  ``Lectures
in Topology'', The University of Michigan Conference of 1940, 
191-221.
\sk1

\ni 47. Nagata, J., Modern Dimension Theory, Wiley (Interscience), N.Y., 1965.
\sk1

\ni 48. Wilder, R.L., Topology of Manifolds, AMS Colloquium
Publications, Vol. 32, 1949.
\sk1

\ni 49. Nagami, K., Dimension Theory, Academic Press,
N.Y., 1970.
\sk1

\ni 50. Bing, R.H. and Floyd, E.E., Coverings with connected
intersections, TAMS 69 (1950) 387-391.
\sk1

\ni 51 Repov{\v s}, D. and {\v S}{\v c}epin, A proof of the Hilbert-Smith
conjecture for actions by Lipschitz maps, Math. Ann. 308 (1997) 361-364.
\sk1

\ni 52 Shchepin, E.V., Hausdorff dimension and the dynamics of
diffeomorphisms, Math Notes 65 (1999) 381-385.
\sk1

\ni 53 Maleshick, Iozke, The Hilbert-Smith conjecture for H{\"o}lder actions,
Usp. Math. Nauk 52 (1997) 173-174.
\sk1

\ni 54 Martin, Gaven J., The Hilbert-Smith conjecture for quasiconformal
actions, Electronic Res. Announcements of the AMS 5 (1999) 66-70.
\vskip 20pt

\cl {The State University of New York at Binghamton}

\end